\documentclass[a4paper,12pt,reqno]{amsart}
\usepackage{amsmath}
\usepackage{amssymb}
\usepackage{amsfonts}
\usepackage{graphicx}
\usepackage{mathtools}
\usepackage[colorlinks]{hyperref}
\renewcommand\eqref[1]{(\ref{#1})} 
\graphicspath{ {images/} }
\title[Hardy type inequalities on Cartan--Hadamard Manifolds]{Critical, stability and higher-order analysis for Hardy type inequalities on Cartan--Hadamard manifolds}
\author[P. Roychowdhury]{Prasun Roychowdhury}
\address{
	Prasun Roychowdhury:
	\endgraf
	Dipartimento di Matematica e Applicazioni
	\endgraf
	Università degli Studi di Milano–Bicocca
	\endgraf
	Via Cozzi 55, 20125 Milano
	\endgraf
	Italy
	\endgraf
	{\it E-mail address} {\rm prasunroychowdhury1994@gmail.com}}
\author[D. Suragan]{Durvudkhan Suragan}
\address{
	Durvudkhan Suragan:
	\endgraf
	Department of Mathematics
	\endgraf
	Nazarbayev University
	\endgraf
	Kazakhstan
	\endgraf
	{\it E-mail address} {\rm durvudkhan.suragan@nu.edu.kz}}
\author[N. Yessirkegenov]{Nurgissa Yessirkegenov}
\address{
  Nurgissa Yessirkegenov:
    \endgraf
  KIMEP University, Almaty, Kazakhstan
  \endgraf
  and
  \endgraf
  Institute of Mathematics and Mathematical Modeling
  \endgraf
Almaty, Kazakhstan
  \endgraf
  {\it E-mail address} {\rm nurgissa.yessirkegenov@gmail.com}
  }
\date{\today}
\theoremstyle{plain}
\newtheorem{theorem}{Theorem}[section]

\newtheorem{lemma}{Lemma}[section]
\newtheorem{corollary}{Corollary}[section]
\newtheorem{remark}{Remark}[section]

\newtheorem{example}{Example}[section]

\numberwithin{equation}{section} \allowdisplaybreaks

\usepackage[text={6in,8.6in},centering]{geometry}
\newcommand{\grad}{\nabla}
\newcommand{\rn}{\mathbb{R}^N}
\newcommand{\dx}{\:{\rm d}x}
\newcommand{\dr}{\:{\rm d}r}
\newcommand{\dv}{\:{\rm d}v_g}

\newcommand{\m}{\mathbb{M}}
\newcommand{\dsn}{\:{\rm d}\sigma}
\newcommand{\ds}{\:{\rm d}s}
\newcommand{\hn}{\mathbb{H}^N}
\newcommand{\sn}{\mathbb{S}^{N-1}}
\newcommand{\ba}{\mathcal{B}_1(o)}
\newcommand{\br}{\mathcal{B}_R(o)}
\newcommand{\ball}[1]{\mathcal{B}_{#1}(o)}
\begin{document}
\begin{abstract}
In this paper, we focus on three main objectives related to Hardy-type inequalities on Cartan--Hadamard manifolds. Firstly, we explore critical Hardy-type inequalities that contain logarithmic terms, highlighting their significance. Secondly, we examine the stability of both critical and subcritical cases of the Hardy inequality. Lastly, we establish two weighted Hardy-type inequalities where singularities appear at the origin as well as the boundary, and we discuss their implications for higher-order operators. Our results improve upon previous findings and also present new higher-order versions of these inequalities as additional outcomes.
\end{abstract}

\maketitle

\tableofcontents

\section{Introduction}\label{intro}
The Cartan--Hadamard manifold  $\m$ is known to admit the following Hardy inequality (see \cite{Carron, YSK}):
\begin{align}\label{int-hardy}
\frac{(N-2)^2}{4}\int_\m\frac{|f(x)|^2}{\rho^2(x)}\dv\leq \int_\m|\nabla_g f(x)|_g^2\dv, \text{ for all }f\in C_0^\infty(\m),
\end{align}
where for $x\in \m$, the nonnegative function $\rho(x)=\text{dist}(x,o)$ denotes the geodesic distance from fixed pole $``o"$, and $\nabla_g$ stands for the Riemannian gradient, $|\cdot|_g$ represents the length of a vector field with respect to the Riemannian metric $g$, and $\dv$ represents volume element in $N\geq 3$ dimensional manifold $(\m,g).$ It is well known that the constant $\frac{(N-2)^2}{4}$ in \eqref{int-hardy} is optimal and never achieved in the underlying Sobolev space, and this suggests improvement is possible. Several authors studied the improvement of this inequality on Riemannian manifolds (we refer to a few of them, e.g., \cite{BGGP, BGR-21, Dambrosio, pinch, GR-23, Kombe1, prc-20, YSK}). Weighted $L^p$-Hardy inequalities have been investigated for both sub-critical and critical cases in \cite{vhn}, and some of our established results here give an immediate improvement of those. 

One of the motivations for this paper is to study the critical behavior of the Hardy inequality on the Cartan--Hadamard manifold. It is natural to expect logarithmic terms in the critical case scenario. A more general form of the critical Hardy inequality on such a manifold was studied in \cite{vhn}. This reads as for all $f\in C_0^\infty(\ba)$, there holds 
\begin{align}\label{int-crt-hardy}
\bigg(\frac{N-1}{N}\bigg)^N\int_{\ba}\frac{|f(x)|^N}{\rho^N(x)\big|\ln\frac{1}{\rho(x)}\big|^{N}}\dv\leq \int_{\ba}\big|\nabla_{g}f(x)\big|_g^N\dv,
\end{align}
with the sharp constant $\big(\frac{N-1}{N}\big)^N$ and $\ba$ is the geodesic ball in $\m$ centered at the pole with unit radius. Here, in this paper, we establish the critical Hardy inequality for the difference between the function and its scaling form. This reads as  for all complex-valued $f\in C_0^\infty(\m)$, there holds 
\begin{align}\label{(1.3)}
\bigg(\frac{N-1}{N}\bigg)^N\int_{\m}\frac{|f(x)-f_1(x)|^N}{\rho^N(x)\big|\ln\frac{1}{\rho(x)}\big|^{N}}\dv\leq \int_{\m}\big|\nabla_{g}f(x)\big|_g^N\dv,  
\end{align}
where $R>0$ is some real number and $f_{1}(x):=f\big(\frac{x}{\rho(x)}\big)$. Moreover, the constant is sharp. In the Euclidean setting, this result was established in \cite{iio}. We prove here more general versions of \eqref{(1.3)} in terms of parameters (see Theorem~\ref{anoth-crit}). Also, we calculate the remainder terms for such inequalities.  

Another type of critical Hardy inequality appears in the literature where the logarithmic term is considered on the other side of \eqref{int-crt-hardy}. This study has been done in a homogeneous group setting in \cite{RSY17_Comptes}, \cite{RSY18_Sobolev_type}, and \cite{RSY18_Tran}. It is natural to address this type of question in a manifold setting. A particular case of our result reads as follows: let $\m$ be a Cartan--Hadamard manifold with curvature bounded above by a negative constant (say $-b$) and $N\geq 2$, then for any complex-valued $f\in C_0^\infty(\m\setminus\{o\})$, there holds
\begin{align}
\int_{\ba^c}\frac{|f(x)|^N}{\rho^N(x)}\dv&\leq \int_{\ba^c}\frac{|f(x)|^N}{\rho^N(x)}\bigg[1+{N(N-1)}D^b(\rho(x)) \, (\ln \rho(x)) \bigg]\dv\nonumber\\&\leq N^N \int_{\ba^c} |\ln \rho(x)|^N \big|\nabla_{g} f(x)\big|_g^N\dv,
\end{align}
where $D^b$ is some nonnegative function related to the density function corresponding to the manifold. 

Moreover, for a special kind of manifold, which satisfies \eqref{condi}, we have the following result: for any complex-valued function $f\in C_0^\infty (\m\setminus \{o\})$, there holds
\begin{align}
\int_{\m}\frac{|f(x)|^N}{\rho^N(x)}\dv\leq N^N \int_{\m} |\ln \rho(x)|^N\big|\nabla_{g}f(x)\big|_g^N\dv,
\end{align}
where the constant $N^N$ is the optimal one. In the case of a flat manifold $\m=\rn$, the required condition is automatically satisfied, see Subsection~\ref{sub-1} for more details. 

Studying the quantitative stability versions of fundamental functional analysis inequalities has become a celebrated topic nowadays. Let us consider the Sobolev inequality
\begin{align}
S_N\bigg(\int_{\rn}|f(x)|^{\frac{2N}{N-2}}\dx\bigg)^{\frac{N-2}{N}}\leq \int_{\rn}|\nabla f(x)|^2\dx \text{ for all }f\in \mathcal{D}_{0}^{1,2}(\rn),
\end{align}
where $\mathcal{D}_{0}^{1,2}(\rn)$ is the completion of $C_0^\infty(\rn)$ function under the norm $||\nabla \cdot||_{L^2(\rn)}$. Also the dimension $N\geq 3$ and $S_N=\frac{N(N-2)}{4}|\sn|^{2/N}$, and $\sn$ denotes the $N$-dimensional unit ball. In \cite{bl}, Brezis and Lieb posed the question of whether it is possible to bound the ``Sobolev deficit" from below in terms of some natural distance (the norm in $\mathcal{D}_{0}^{1,2}(\rn)$) from the set of optimizers i.e., do there exist $\kappa,\, \alpha>0$ such that 
\begin{align*}
\kappa \leq \frac{S_N\bigg(\int_{\rn}|f(x)|^{\frac{2N}{N-2}}\dx\bigg)^{\frac{N-2}{N}}- \int_{\rn}|\nabla f(x)|^2\dx}{\text{dist}(
f,\mathcal{M})^\alpha },
\end{align*}
where $f\in \mathcal{D}_{0}^{1,2}(\rn)$ and the optimizer set is defined as
\begin{align*}
\mathcal{M}&:=\{f_{a,b,c}\, : (a,b,c)\in (0,+\infty)\times \rn \times \mathbb{R}\},\\& \text{where } f_{a,b,c}(x)=c\bigg(\frac{2}{1+|\frac{x-b}{a}|^2}\bigg)^{\frac{N-2}{2}}?
\end{align*}
This question was answered by Bianchi and Egnell in \cite{be}. It was shown by Lions \cite{lion} that, if the Sobolev deficit is small for some function $f$, then $f$ has to be close to the set of Sobolev optimizers. The closeness is measured with respect to the norm in $\mathcal{D}_{0}^{1,2}(\rn)$, which is the strongest possible sense. The Bianchi--Egnell result makes the qualitative result of Lions more quantitative. In recent years, stability analysis has gained significant attention. For a list of related works, we refer to \cite{cfmp, cgm, deffl, fg, fn, f-duke, rlf, lw, ne} and the references therein.

A similar question can be asked for the Hardy inequality. But the main challenge here is that one does not have an extremizer of \eqref{int-hardy}. However, it is known that the sharpness of the Hardy constant $\frac{(N-2)^2}{4}$ can be shown by using the approximation of the function $\rho^{-\frac{(N-2)}{2}}(x)$. Exploiting some scale versions of this function, we analyse the stability of the Hardy inequality. We define for the function  $f\in C_0^\infty(\m\setminus\{o\})$ and for a fixed real number $R>0$, the following distance function
\begin{align*}
d_H(f,R):= \bigg(\int_{\m}\frac{\big|f(x)-R^{\frac{N-2}{2}}f(R\frac{x}{\rho(x)})\rho^{-\frac{N-2}{2}}(x)\big|^2}{|\ln \frac{R}{\rho(x)}|^2\rho^{2}(x)}\dv\bigg)^{\frac{1}{2}}.
\end{align*}
Now our result states that for $N\geq 3$ and complex-valued function $f\in C_0^\infty(\m\setminus\{o\})$ there holds
\begin{align}
\frac{1}{4} \sup_{R>0} d_H(f,R)^2 \leq \int_{\m}|\nabla_{g}f(x)|_g^2\dv- \frac{(N-2)^2}{4}\int_{\m}\frac{|f(x)|^2}{\rho^{2}(x)}\dv.
\end{align}
We prove the results for more general weighted $L^p$-Hardy inequalities for critical and sub-critical cases (see Subsection~\ref{sub-2} for details).

In the final section of our note, we establish an improvement of both classical and geometric Hardy inequalities. This type of result was studied in \cite{sano} for the Euclidean case and recently extended for the homogeneous group in \cite{SY-arxiv}. By tackling the density function here, we study the result for the Cartan--Hadamard manifold setting. We establish the following result:  let $2\leq b<\infty$, $a<N$, $c>0$, and $a\leq N - (b-1)c$. Then for any complex-valued function $f\in C_0^\infty (\ba)$ there holds
\begin{align}\label{int-two-hardy}
\bigg(\frac{(b-1)}{2}c\bigg)^2\int_{\ba} \frac{|f(x)|^{2}}{\rho^{a}(x)(1-\rho^{c}(x))^{b}} \dv\leq \int_{\ba} \frac{|\nabla_{g} f(x)|_g^2}{\rho^{a-2}(x)(1-\rho^{c}(x))^{b-2}} \dv,
\end{align}   
where the constant $\big(\frac{(b-1)}{2}c\big)^2$ in this inequality is sharp. As an immediate implication, one can see by selecting $a=b=2$ and  $c=N-2$ in \eqref{int-two-hardy}, for any $f\in C_0^\infty (\ba)$ with $N\geq 3$, there holds
\begin{align}
\frac{(N-2)^2}{4}\int_{\ba} \frac{|f(x)|^{2}}{\rho^{2}(x)} \dv&\leq  \frac{(N-2)^2}{4}\int_{\ba} \frac{|f(x)|^{2}}{\rho^{2}(x)(1-\rho^{N-2}(x))^{b}} \dv \nonumber \\&\leq \int_{\ba} |\nabla_{g} f(x)|_g^2\dv.
\end{align}
This gives an improvement of the classical Hardy inequality \eqref{int-hardy}. Now let us consider the distance function $\text{dist}(x,\partial\ba)=1-\rho(x)$ and substitute $b=a=2$ and $c=1$ in \eqref{int-two-hardy}, then we have the following improvement of the geometric Hardy inequality. This reads as follows: for any $f\in C_0^\infty (\ba)$ with $N\geq 3$ there holds
\begin{align}\label{geo}
\frac{1}{4}\int_{\ba}\frac{|f(x)|^2}{\text{dist}^2(x,\partial\ba)}\dv&\leq \frac{1}{4}\int_{\ba}\frac{|f(x)|^2}{\rho^2(x)\text{dist}^2(x,\partial\ba)}\dv\nonumber \\&\leq \int_{\ba} |\nabla_{g} f(x)|_g^2 \dv.
\end{align}
Moreover, the constant $\frac{1}{4}$ in \eqref{geo} is sharp. In general, the best constant depends on the geometry of the domain. For more details, we refer \cite{bar-1, bar-2, bm, LLZ}, to name a few. Moreover, our result holds for any $p\in (
1,\infty)$ (for the details see Subsection~\ref{sub-3}). 

As an application, we obtain the Caffarelli--Kohn--Nirenberg inequality involving two singularities. Note that Rellich inequalities are extensions of Hardy inequalities in higher-order Sobolev spaces. In that spirit, we develop new inequalities for the higher-order operator that improve previously known results.  Let us mention the improvement of the Rellich inequality here (see the discussions in Subsection~\ref{sub-4} for details).  Let $(\m,g,\psi)$ be an $N$-dimensional Riemannian model manifold with $N\geq 5$ and $-\frac{\psi^{\prime\prime}(r)}{\psi(r)}\leq -1$ (for detailed definitions see Section~\ref{sec-2}). Then for any $f\in C_0^\infty (\ba)$ there holds
\begin{align}
\frac{N^2(N-4)^2}{16}\int_{\ba} \frac{|f(x)|^{2}}{\rho^{4}(x)} \dv &\leq  \frac{N^2(N-4)^2}{16}\int_{\ba} \frac{|f(x)|^{2}}{\rho^{4}(x)(1-\rho^{N-4}(x))^{2}} \dv \nonumber\\& \leq  \int_{\ba} |\Delta_{g} f(x)|^2 \dv,
\end{align}
where $\Delta_g$ represents the Laplace-Beltrami operator and the constant $\frac{N^2(N-4)^2}{16}$ is the best possible.

The structure of this paper is outlined as follows. Section~\ref{sec-2} delves into the fundamental aspects of Cartan--Hadamard manifolds. Following this, Section~\ref{sec-3} delineates our main results, which are further subdivided into four distinct subsections. These cover a range of topics, including the critical Hardy inequalities, stability analysis, the Hardy inequality with a dual singularity, and the advancement of higher-order results that enhance the established Rellich inequality on manifolds.

In Section~\ref{pfs-1}, we validate the critical Hardy inequality using logarithmic terms. Section~\ref{pfs-2} is dedicated to proving an alternate form of the critical Hardy inequality, where logarithmic terms are positioned on the right-hand side of the inequality. Section~\ref{pfs-3} explores the stability aspect of the Hardy inequalities.

Moving on, Section~\ref{pfs-4} focuses on establishing two weighted Hardy inequalities that incorporate singularities both at the pole and at the boundary. Section~\ref{pfs-5} discusses the higher-order results, emphasizing the interesting role of basic induction. We conclude the paper in Section~\ref{sec-open} by proposing some open questions.

\medspace

\section{Preliminaries}\label{sec-2}
This section will briefly recall some of the known facts concerning Cartan--Hadamard manifolds. Let $(\mathbb{M},g)$ be an $N$-dimensional complete Riemannian manifold. In terms of the local coordinate system $\{x^i\}_{i=1}^{N},$ one can write
\begin{align*}
g=\sum g_{ij}{\rm d}x^i\dx^j.
\end{align*}
The Laplace-Beltrami operator $\Delta_g$ concerning the metric $g$ is defined as follows
\begin{align*}
\Delta_g:=\sum \frac{1}{\sqrt{\text{det }(g_{ij})}}\frac{\partial}{\partial x^i}\bigg(\sqrt{\text{det }g_{ij}}\:g^{ij}\frac{\partial}{\partial x^j}\bigg),
\end{align*}
where $(g^{ij})=(g_{ij})^{-1}$. Also, we denote $\nabla_g$ as the Riemannian gradient, and we have\begin{align*}
\langle \nabla_g f , \nabla_g h \rangle_g= \sum g^{ij}\frac{\partial f}{\partial x^i}\frac{\partial h}{\partial x^j}.
\end{align*}
For simplicity, we shall use the notation $|\nabla_g f|_g = \sqrt{g(\nabla_g f,\nabla_g f)}$.

We say a Riemannian manifold $(\mathbb{M},g)$ is a Cartan--Hadamard manifold if it is complete, simply connected, and admits a non-positive sectional curvature, i.e., the sectional curvature $\text{K}_{\m}\leq 0$ along each plane section at each point of $\m$. Let us fix a point $o\in M$ called the pole or the origin of the manifold $\m$. Let $\rho(x)=\text{dist}(x,o)$ be the geodesic distance between $x\in \m$ and $o$. Note that, $\rho(x)$ is smooth on $\m\setminus\{o\}$ and for all $x\in \m\setminus\{o\}$ satisfies $|\nabla_g \rho(x)|_g=1$. For the Cartan--Hadamard manifold $(\m,g)$, the exponential map $\text{Exp}_o\::\: \text{T}_o\m\rightarrow \m$ is a diffeomorphism, where $\text{T}_o\m$ is the tangent space to $\m$ at the pole $o$. 

For any function $f\in C^1(\m)$, we define the radial derivation $\partial_\rho=\frac{\partial}{\partial \rho}$ along the geodesic curve starting from pole $o$ by 
\begin{align*}
\partial_\rho f(x) = \frac{{\rm d}(f\circ \text{Exp}_o)}{{\rm d}r}\big(\text{Exp}_o^{-1}(x)\big),
\end{align*}
where $\frac{{\rm d}}{{\rm d}r}$ denotes the radial derivation on $T_o\m$, i.e.,
\begin{align*}
\frac{{\rm d}}{{\rm d}r}F(\sigma)=\bigg\langle \frac{\sigma}{|\sigma|}\:,\:\nabla F(\sigma)\bigg\rangle, \quad \quad \sigma\in T_o\m\setminus\{0\}.
\end{align*}
We note that by Gauss's lemma (see, e.g., \cite{rim-2}), we have that 
\begin{align}\label{gauss}
|\partial_\rho f |\leq |\nabla_{g} f|_g.
\end{align}

For any $\delta>0$, we denote by $\mathcal{B}_\delta(o):=\{x\in \m\,:\, \text{dist}(x,o)<\delta \}$ the geodesic ball in $\m$ with center at fixed pole $o$ and radius $\delta$. For any function $f\in L^1(\m)$, the polar coordinate decomposition can be written as follows
\begin{align*}	\int_{\m}f(x)\dv=\int_{\sn}\int_{0}^{\infty}f(\text{Exp}_o(t\sigma))\:J(t,\sigma)t^{N-1}{\rm d}t\dsn,
\end{align*}
where $\dsn$ denotes the measure on the unit sphere of the tangent space $\text{T}_o\m$ and $J(t,\sigma)$  is the density function for the manifold $(\m,g)$. For the details, we refer to \cite{rim-2,rim-1}. Thus, we have the expression of the Jacobian term in the form 
\begin{align*}
\lambda(t,\sigma)=J(t,\sigma)\:t^{N-1}.
\end{align*}
Also, see \cite[p.~166-169]{rim-2}, it is known that $J(t,\sigma)$ is increasing and there holds
\begin{align}\label{imp-den}
J(t,\sigma)=1+O(t^2) \text{ for } t\rightarrow 0^+.
\end{align}

Recall the radial Laplace operator on $N$-dimensional Cartan--Hadamard manifold $(\m,g)$ by 
\begin{align*}
\Delta_{\rho,g}f(x)=\partial^2_\rho f(x)+\bigg(\frac{N-1}{\rho(x)}+\frac{J_\rho(\rho(x),\sigma)}{J(\rho(x),\sigma)}\bigg)\partial_\rho f(x), \quad x=\text{Exp}_o(\rho \sigma),
\end{align*}
with $f\in C^2(\m)$, where the lower subscript $\rho$ denotes the derivative with respect to the radial component $\rho$. 

Let us consider Cartan--Hadamard manifolds that admit constant curvature, i.e., $\text{K}_{\m} \equiv -b$ for some $b\geq 0$. Then it is clear that $J(t,\sigma)=J^b(t)$, where 
\begin{equation}
J^b(t)=\;
\begin{dcases}
		1 & \text{ if } b=0, \\
		\bigg(\frac{\sinh (\sqrt{b}t)}{\sqrt{b}t}\bigg)^{N-1}  & \text{ if } b>0. \\
	\end{dcases}
\end{equation}

Now, let us focus on the Cartan--Hadamard manifolds whose sectional curvature is $K_\m\leq -b$ for some $b\geq 0$.  First, for $b\geq 0$, we consider the function $C^b :(0,\infty)\rightarrow \mathbb{R}$ defined by
\begin{equation*}
	C^b(t)=\;
	\begin{dcases}
		\frac{1}{t} & \text{ if } b=0, \\
		\sqrt{b} \coth (\sqrt{b} t)  & \text{ if } b>0, \\
	\end{dcases}
\end{equation*}
and the function $D^b :[0,\infty)\rightarrow \mathbb{R}$ defined by
\begin{equation*}
	D^b(t)=\;
	\begin{dcases}
		0 & \text{ if } b=0, \\
		t C^b(t)-1  & \text{ if } b>0. \\
	\end{dcases}
\end{equation*}
Then, the Bishop--Gunther comparison theorem  \cite[p.~172]{rim-2} says that
\begin{align}\label{density_comp}
 (\ln J(t,\sigma))_t=\frac{J_t(t,\sigma)}{J(t,\sigma)}\geq \frac{J^b_r(t)}{J^b(t)}   = \frac{N-1}{t} D^b(t), \text{ for } t>0.
\end{align}
$D^b(t)$ is a nonnegative function and therefore for any Cartan--Hadamard manifolds, we have 
\begin{align}\label{den-der}
\frac{J_t(t,\sigma)}{J(t,\sigma)}\geq 0.
\end{align}

To conclude this section, we examine a few explicit examples of Cartan--Hadamard manifolds.
\begin{example}
	If $K_{\m}\equiv 0$, then $\m=\rn$ (Euclidean space) and $\lambda(t,\sigma)=t^{N-1}$.
\end{example}

\begin{example}
	If $K_{\m}\equiv -1$, then $\m=\hn$ (Hyperbolic space) and the Jacobian is given by $\lambda(t,\sigma)=(\sinh t)^{N-1}$, and so the density function will be $J(t,\sigma)=\big(\frac{\sinh t}{t}\big)^{N-1}$.
\end{example} 

\begin{example}
An $N$-dimensional Riemannian symmetric model manifold $(\m,g,\psi)$ is a Riemannian manifold whose metric $g$ is represented in spherical coordinates by
\begin{equation*}
	{\rm d}s^2 = {\rm d}t^2 + \psi^2(t) \, {\rm d}\sigma^2,
\end{equation*}
where the coordinate $t$, denotes the Riemannian distance from the fixed pole $o$, and ${\rm d}\sigma^2$ is the metric on $N$-dimensional sphere $\mathbb{S}^{N-1}$ and the non-negative function $\psi$ satisfies
\begin{align*}
	\psi\in C^\infty([0,+\infty)), \,  \psi > 0 \text{ on } (0,+\infty),  \psi'(0) =1 \text{ and } \psi^{(2k)}(0)= 0 \text{ for all } k\geq 0\,,
\end{align*}
where prime denotes the derivative w.r.t. $t$. It is known that for this model manifold $(\m,g,\psi)$, the density function $J(r,\sigma)=\big(\frac{\psi(t)}{t}\big)^{N-1}$, and hence we have the Jacobian $\lambda(t,\sigma)=(\psi(t))^{N-1}$. Moreover, the sectional curvature is given by $\text{K}_{\m}=-\frac{\psi^{\prime\prime}(t)}{\psi(t)}$. Therefore, $(\m,g,\psi)$ is a Cartan--Hadamard manifold when we have $\psi^{\prime\prime}(t)\geq 0$ on $(0,\infty)$. For example, when $\psi(t)=t$ and $\psi(t)=(\sinh t)$, they coincide with the known Cartan--Hadamard manifold as Euclidean space and Hyperbolic space, respectively.    
\end{example}

\medspace

\section{Statement of main results}\label{sec-3}
This section is devoted to the statements of the main results. We consider Cartan--Hadamard manifolds $(\m,g)$ with dimension $N\geq 2$. Note that from now onwards $C_0^\infty(\Omega)$ and $L^p(\Omega)$ mean compactly supported smooth functions and $L^p$ integrable function space on $\Omega$, respectively,  where $\Omega\subset \m$ is a domain with smooth boundary or a domain without boundary. We denote by $W^{1,p}_0(\Omega)$ the Sobolev space consisting of functions in $L^p(\Omega)$ whose first-order weak derivatives also belong to $L^p(\Omega)$ and whose trace vanishes on the boundary.

\subsection{Critical Hardy inequalities}\label{sub-1} First, we present logarithmic type critical Hardy inequalities.

\begin{theorem}\label{anoth-crit-th}
    Let $1<\gamma<\infty$ and $\max\{1,\gamma-1\}<p<\infty$. Then for all complex-valued functions $f\in C_0^\infty(\m)$ and all $R>0$ we have the inequality
    \begin{align}\label{anoth-crit}
        \frac{(\gamma-1)}{p}\bigg|\bigg|\frac{f(x)-f_R(x)}{\rho^{\frac{N}{p}}(x)\big(\ln\frac{R}{\rho(x)}\big)^{\frac{\gamma}{p}}}\bigg|\bigg|_{L^p(\m)}\leq \bigg|\bigg|\rho^{\frac{p-N}{p}}(x)\bigg(\ln \frac{R}{\rho(x)}\bigg)^{\frac{p-\gamma}{p}}\partial_\rho f(x)\bigg|\bigg|_{L^p(\m)},
    \end{align}
    and applying \eqref{gauss} we have
    \begin{align}\label{ful-anoth-crit}
        \frac{(\gamma-1)}{p}\bigg|\bigg|\frac{f(x)-f_R(x)}{\rho^{\frac{N}{p}}(x)\big(\ln\frac{R}{\rho(x)}\big)^{\frac{\gamma}{p}}}\bigg|\bigg|_{L^p(\m)}\leq \bigg|\bigg|\rho^{\frac{p-N}{p}}(x)\bigg(\ln \frac{R}{\rho(x)}\bigg)^{\frac{p-\gamma}{p}}|\nabla_{g} f(x)|_g\bigg|\bigg|_{L^p(\m)},
    \end{align}
    where $f_R(x):=f\big(R\frac{x}{\rho(x)}\big)$ and the constant $\frac{(\gamma-1)}{p}$ is sharp.
\end{theorem}

The inequality \eqref{anoth-crit} was studied in \cite[Theorem~2.2.7]{rs-book} when the underlying space was a homogeneous group $\mathbb{G}$. It is worth mentioning that a higher-order version of \eqref{anoth-crit} on the homogeneous group $\mathbb{G}$ was studied in \cite[Theorem~5.2]{vhn-2}. We refer the reader to the same article for a comprehensive review of many Hardy--Rellich type inequalities. If we observe the steps of the proof of Theorem~\ref{anoth-crit-th} (refer to Section~\ref{pfs-1}) and apply the H\"older inequality, we deduce an improvement of \eqref{anoth-crit} with a remainder term and involving the density function of the related manifold. This reads as follows:
\begin{corollary}\label{anoth-crit-cor-rem}
    Let $1<\gamma<\infty$ and $\max\{1,\gamma-1\}<p<\infty$. Then for all complex-valued functions $f\in C_0^\infty(\m)$ and all $R>0$ we have the inequality
    \begin{align*}
       \frac{\bigg|\bigg|\frac{(f(x)-f_R(x))(\partial_\rho\ln(J(\rho(x),\sigma))^{1/p}}{\rho^{\frac{N-1}{p}}(x)\big(\ln \frac{R}{\rho(x)}\big)^{\frac{\gamma-1}{p}}}\bigg|\bigg|^p_{L^p(\m)}}{p\bigg|\bigg|\frac{f(x)-f_R(x)}{\rho^{\frac{N}{p}}(x)\big(\ln\frac{R}{\rho(x)}\big)^{\frac{\gamma}{p}}}\bigg|\bigg|_{L^p(\m)}^{(p-1)}}&+\frac{(\gamma-1)}{p}\bigg|\bigg|\frac{f(x)-f_R(x)}{\rho^{\frac{N}{p}}(x)\big(\ln\frac{R}{\rho(x)}\big)^{\frac{\gamma}{p}}}\bigg|\bigg|_{L^p(\m)}\\&\leq \bigg|\bigg|\rho^{\frac{p-N}{p}}(x)\bigg(\ln \frac{R}{\rho(x)}\bigg)^{\frac{p-\gamma}{p}}\partial_\rho f(x)\bigg|\bigg|_{L^p(\m)},
    \end{align*}
    and applying \eqref{gauss} we have
\begin{align*}
\frac{\bigg|\bigg|\frac{(f(x)-f_R(x))(\partial_\rho\ln(J(\rho(x),\sigma))^{1/p}}{\rho^{\frac{N-1}{p}}(x)\big(\ln \frac{R}{\rho(x)}\big)^{\frac{\gamma-1}{p}}}\bigg|\bigg|^p_{L^p(\m)}}{p\bigg|\bigg|\frac{f(x)-f_R(x)}{\rho^{\frac{N}{p}}(x)\big(\ln\frac{R}{\rho(x)}\big)^{\frac{\gamma}{p}}}\bigg|\bigg|_{L^p(\m)}^{(p-1)}}&+\frac{(\gamma-1)}{p}\bigg|\bigg|\frac{f(x)-f_R(x)}{\rho^{\frac{N}{p}}(x)\big(\ln\frac{R}{\rho(x)}\big)^{\frac{\gamma}{p}}}\bigg|\bigg|_{L^p(\m)}\\&\leq \bigg|\bigg|\rho^{\frac{p-N}{p}}(x)\bigg(\ln \frac{R}{\rho(x)}\bigg)^{\frac{p-\gamma}{p}}|\nabla_g f(x)|_g\bigg|\bigg|_{L^p(\m)},
\end{align*}
    where $f_R(x):=f\big(R\frac{x}{\rho(x)}\big)$ and the constant $\frac{(\gamma-1)}{p}$ is sharp.
\end{corollary}

If we substitute $\gamma=p$ in \eqref{anoth-crit}, we have the following result. This result will be used for further context in this note.
\begin{remark}
    Let $1<p<\infty$. Then for all complex-valued functions $f\in C_0^\infty(\m)$ and all $R>0$ we have the inequality
    \begin{align}\label{anoth-crit-cor}
        \bigg(\frac{p-1}{p}\bigg)^p\int_{\m}\frac{|f(x)-f_R(x)|^p}{\rho^{N}(x)\big|\ln\frac{R}{\rho(x)}\big|^{p}}\dv\leq \int_{\m}\rho^{p-N}(x)\big|\partial_\rho f(x)\big|^p\dv,
    \end{align}
      and applying \eqref{gauss} we have
      \begin{align}
        \bigg(\frac{p-1}{p}\bigg)^p\int_{\m}\frac{|f(x)-f_R(x)|^p}{\rho^{N}(x)\big|\ln\frac{R}{\rho(x)}\big|^{p}}\dv\leq \int_{\m}\rho^{p-N}(x)\big|\nabla_{g}f(x)\big|_g^p\dv,
    \end{align}
    where $f_R(x):=f\big(R\frac{x}{\rho(x)}\big)$ and the constant $\big(\frac{p-1}{p}\big)^p$ is sharp.
\end{remark}

Now let us focus on critical Hardy inequalities. On the Cartan--Hadamard manifolds, these types of inequalities have been studied in \cite[Theorem~3.3]{vhn}. Note that the logarithmic term has not been considered with the gradient term there. Here, we have first made an effort to establish this type of result. We will use the following functional from \cite[Proposition~1.1]{iio}: 
\begin{equation}\label{Rp_1}
R_p(\zeta,\eta):=\begin{cases} {}\left(\frac{1}{p}|\eta|^p+\frac{p-1}{p}|\zeta|^p-|\zeta|^{p-2} \operatorname{Re}(\zeta \bar{\eta})\right)/|\zeta-\eta|^{2} \;\text{\;if}\;\;\zeta \neq \eta,\\
\frac{p-1}{2}|\zeta|^{p-2}\;\text{\;if}\;\;\zeta=\eta\end{cases}
\end{equation}
for any $1<p<\infty$ and any complex-valued $\zeta,\eta$. Then we have
\begin{theorem}\label{new_critical_hardy}
Let $1<p<\infty$. Then there holds the identity for any complex-valued function $f\in C_0^\infty (\m\setminus \{o\})$
\begin{multline}\label{new_critical_hardy_1}
    \int_{\m}\frac{|f(x)|^p}{\rho^N(x)}\dv=p^p \int_{\m} \frac{|\ln \rho(x)|^p}{\rho^{N-p}(x)}|\partial_\rho f(x)|^p\dv \\ - p\int_{\m} |f(x)|^p \frac{J_\rho(\rho(x),\sigma)}{J(\rho(x),\sigma)}\, (\ln \rho(x)) \rho^{1-N}(x)\dv \\-p\int_{\m} R_p\bigg(\rho^{-\frac{N}{p}}(x)f(x),-p\frac{(\ln \rho(x))}{\rho^{\frac{(N-p)}{p}}(x)}\partial_\rho f(x)\bigg)\left|\rho^{-\frac{N}{p}}(x)f(x)+p\frac{(\ln \rho(x))}{\rho^{\frac{(N-p)}{p}}(x)}\partial_\rho f(x)\right|^{2}\dv,
\end{multline}
where $R_{p}$ is defined in \eqref{Rp_1}.
Moreover, we have the following inequality for any complex-valued function $f\in C_0^\infty (\m\setminus \{o\})$
\begin{multline}\label{new_critical_hardy_2}
  \left(\int_{\m}\frac{|f(x)|^p}{\rho^N(x)}\dv\right)^{\frac{1}{p}} \leq p\left(\int_{\m}|\partial_{\rho}f(x)|^{p}\frac{|\ln \rho(x)|^{p}}{\rho^{N-p}(x)}\dv\right)^{\frac{1}{p}}\\- \left(\int_{\m} |f(x)|^p \frac{J_\rho(\rho(x),\sigma)}{J(\rho(x),\sigma)}\, (\ln \rho(x)) \rho^{1-N}(x)\dv\right) \left(\int_{\m}\frac{|f(x)|^p}{\rho^N(x)}\dv\right)^{\frac{1-p}{p}},
\end{multline}
and applying \eqref{gauss}, we have
\begin{multline}\label{new_critical_hardy_3}
   \left(\int_{\m}\frac{|f(x)|^p}{\rho^N(x)}\dv\right)^{\frac{1}{p}} \leq p\left(\int_{\m}|\nabla_{g}f(x)|_g^{p}\frac{|\ln \rho(x)|^{p}}{\rho^{N-p}(x)}\dv\right)^{\frac{1}{p}}\\- \left(\int_{\m} |f(x)|^p \frac{J_\rho(\rho(x),\sigma)}{J(\rho(x),\sigma)}\, (\ln \rho(x)) \rho^{1-N}(x)\dv\right) \left(\int_{\m}\frac{|f(x)|^p}{\rho^N(x)}\dv\right)^{\frac{1-p}{p}}.
\end{multline}
\end{theorem}
\begin{remark}
 Note that by \cite[Proposition~1.1]{iio} for any $2\leq p < \infty$ and all complex-valued $\zeta$ and $\eta$ we have $ R_p(\zeta,\eta)\geq 0$. In general, for all $1<p<\infty$ we have
$$
R_p(\zeta,\eta)=(p-1)\int_{0}^{1}|t\zeta +(1-t)\eta|^{p-2}t \,{\rm d}t \geq 0
$$ 
for any real-valued $\zeta$ and $\eta$ (see \cite[Proposition~1]{iio-17}).
\end{remark}
By this remark, one can observe that we can drop the last integral in \eqref{new_critical_hardy_1} in the following two cases: when $2\leq p<\infty$ for any complex-valued $f\in C_0^\infty (\m\setminus \{o\})$ or when $1<p<\infty$ but for any real-valued $f\in C_0^\infty (\m\setminus \{o\})$. Thus, in these cases, dropping the last integral in \eqref{new_critical_hardy_1} and using \eqref{gauss}, we obtain
\begin{equation}\label{complex_p2_thm3.2_1}
\int_{\m}\frac{|f(x)|^p}{\rho^N(x)}\bigg[1+p\partial_\rho(\ln(J(\rho(x),\sigma))\, (\ln \rho(x)) \rho(x)\bigg]\dv\leq p^p \int_{\m} \frac{|\ln \rho(x)|^p}{\rho^{N-p}(x)}\big|\partial_\rho f(x)\big|^p\dv,
\end{equation}
and 
\begin{equation}\label{complex_p2_thm3.2_2}
\int_{\m}\frac{|f(x)|^p}{\rho^N(x)}\bigg[1+p\partial_\rho(\ln(J(\rho(x),\sigma))\, (\ln \rho(x)) \rho(x)\bigg]\dv\leq p^p \int_{\m} \frac{|\ln \rho(x)|^p}{\rho^{N-p}(x)}\big|\nabla_g f(x)\big|_g^p\dv.
\end{equation}
For the Euclidean space $\m=\rn$, the density function is $J(t,\sigma)=1$, and using this in \eqref{new_critical_hardy_2} of Theorem~\ref{new_critical_hardy} we obtain the critical Hardy inequality which can be compared with \cite[Theorem~2.2.4]{rs-book} for the $\mathbb{G}=(\rn, +)$ case. Also, note that this result holds when the density function is only dependent on the spherical part. Precisely, those Cartan--Hadamard manifolds whose density function satisfies the following condition in polar coordinate representation:
\begin{align}\label{condi}
    J(t,\sigma)=J(\sigma)  \text{ for all } t\in (0,\infty) \text{ and } \sigma\in \sn.
\end{align}
Clearly, $\rn$ is an example of such a manifold. We have the following corollary which follows from \eqref{new_critical_hardy_2} and \eqref{new_critical_hardy_3}:
\begin{corollary}\label{new-critical-hardy-rem}
 Let $1<p<\infty$ and $(\m,g)$ be a Cartan--Hadamard manifold which satisfy \eqref{condi}. Then for any complex-valued function $f\in C_0^\infty (\m\setminus \{o\})$ there holds
\begin{equation}\label{complex_p2_thm3.2_3}
\int_{\m}\frac{|f(x)|^p}{\rho^N(x)}\dv\leq p^p \int_{\m} \frac{|\ln \rho(x)|^p}{\rho^{N-p}(x)}\big|\partial_\rho f(x)\big|^p\dv,
\end{equation}
and applying \eqref{gauss} we have
\begin{equation}\label{complex_p2_thm3.2_4}
\int_{\m}\frac{|f(x)|^p}{\rho^N(x)}\dv\leq p^p \int_{\m} \frac{|\ln \rho(x)|^p}{\rho^{N-p}(x)}\big|\nabla_g f(x)\big|_g^p\dv.
\end{equation} Moreover, the constant $p^p$ is sharp.
\end{corollary}

We are exploiting the simple fact that $(\text{ln} \, t)\geq 0$ for $t\in(1,\infty)$, $D_b\geq 0$ and \eqref{density_comp} in \eqref{complex_p2_thm3.2_1} and \eqref{complex_p2_thm3.2_2}, we have the following immediate consequence on the exterior of the unit ball centered at pole $o$.
\begin{remark}
    Suppose we consider the Cartan--Hadamard manifold $(\m,g)$ with the sectional curvature $K_\m\leq -b$, $b\geq 0$. Then for any $2\leq p<\infty$ and any complex-valued function $f\in  C_0^\infty (\ba^c)$ there holds
\begin{multline}\label{complex_p2_thm3.2_5}
    \int_{\ba^c}\frac{|f(x)|^p}{\rho^N(x)}\dv\leq\int_{\ba^c}\frac{|f(x)|^p}{\rho^N(x)}\bigg[1+{p(N-1)}D^b(\rho(x)) \, (\ln \rho(x)) \bigg]\dv\\\leq p^p \int_{\ba^c} \frac{|\ln \rho(x)|^p}{\rho^{N-p}(x)}\big|\partial_\rho f(x)\big|^p\dv,
\end{multline}
and applying \eqref{gauss} we have
\begin{multline}\label{complex_p2_thm3.2_6}
    \int_{\ba^c}\frac{|f(x)|^p}{\rho^N(x)}\dv\leq\int_{\ba^c}\frac{|f(x)|^p}{\rho^N(x)}\bigg[1+{p(N-1)}D^b(\rho(x)) \, (\ln \rho(x)) \bigg]\dv\\\leq p^p \int_{\ba^c} \frac{|\ln \rho(x)|^p}{\rho^{N-p}(x)}\big|\nabla_g  f(x)\big|_g^p\dv.
\end{multline} Moreover, the inequalities \eqref{complex_p2_thm3.2_5} and \eqref{complex_p2_thm3.2_6} hold for any $1<p<\infty$ and any real-valued function $f\in  C_0^\infty (\ba^c)$.
\end{remark}

\subsection{Stability of Hardy inequality}\label{sub-2} The concept of stability in inequality, often referred to as estimating remainder terms by considering the function's proximity to the set of extremisers, is a fundamental idea. In this discussion, we explore the stability of the Hardy inequality on the Cartan--Hadamard manifold, addressing both critical and sub-critical scenarios. Specifically, we revisit the sub-critical weighted $L^p$-Hardy inequality, as introduced in \cite[Theorem~3.1]{vhn}, within the context of the Cartan--Hadamard manifold. This reads as follows: let $(\m,g)$ be an N-dimensional Cartan--Hadamard manifold. Suppose that $N\geq 2$, $p\in (1,N)$ and $p+\beta<N$. There holds for any $f\in C_0^\infty(\m\setminus\{o\})$
	\begin{align*}
		\bigg(\frac{N-p-\beta}{p}\bigg)^p\int_{\m}\frac{|f(x)|^p}{\rho^{p+\beta}(x)}\dv \leq  \int_{\m}\frac{|\partial_\rho f(x)|^p}{\rho^\beta(x)}\dv,
	\end{align*}
 where the constant $\big(\frac{N-p-\beta}{p}\big)^p$ is the best possible. For the function  $f\in C_0^\infty(\m\setminus\{o\})$ and for a fixed real number $R>0$ we define the following distance function
\begin{align*}
	d_H(f,R):= \bigg(\int_{\m}\frac{\big|f(x)-R^{\frac{N-p-\beta}{p}}f(R\frac{x}{\rho(x)})\rho^{-\frac{N-p-\beta}{p}}(x)\big|^p}{|\ln \frac{R}{\rho(x)}|^p \rho^{p+\beta}(x)}\dv\bigg)^{\frac{1}{p}}.
\end{align*}
Now we are ready to present the following stability result in the sub-critical setting:
\begin{theorem}\label{stab-th}
Let $2\leq p < N$ and $p+\beta<N$. Then for all complex-valued functions $f\in C_0^\infty(\m\setminus\{o\})$, we have 
	\begin{align*}
		c_p\bigg(\frac{p-1}{p}\bigg)^p \sup_{R>0} d_H(f,R)^p \leq \int_{\m}\frac{|\partial_\rho f(x)|^p}{\rho^\beta(x)}\dv- \bigg(\frac{N-p-\beta}{p}\bigg)^p\int_{\m}\frac{|f(x)|^p}{\rho^{p+\beta}(x)}\dv,
	\end{align*}
 and applying \eqref{gauss} we have 
 \begin{align*}
		c_p\bigg(\frac{p-1}{p}\bigg)^p \sup_{R>0} d_H(f,R)^p \leq \int_{\m}\frac{|\nabla_g f(x)|_g^p}{\rho^\beta(x)}\dv- \bigg(\frac{N-p-\beta}{p}\bigg)^p\int_{\m}\frac{|f(x)|^p}{\rho^{p+\beta}(x)}\dv,
\end{align*}
where $c_p$ is mentioned in Lemma~\ref{cpl}.
\end{theorem}

To discuss the stability of the $L^p$-critical Hardy inequality, let us recall the inequality \cite[Theorem~3.3]{vhn}: let $(\m,g)$ be an N-dimensional Cartan--Hadamard manifold with $N\geq 2$. Then for any $p\in(1,\infty)$, there holds 
\begin{align*}
		\bigg(\frac{p-1}{p}\bigg)^p\int_{\ba}\frac{|f(x)|^p}{\rho^{N}(x)(\ln \frac{1}{\rho(x)})^p}\dv \leq \int_{\ba}\frac{|\partial_\rho f(x)|^p}{\rho^{N-p}(x)}\dv,
\end{align*}
where the constant $\big(\frac{p-1}{p}\big)^p$ is sharp. Then we define for a complex-valued function  $f\in C_0^\infty(\ba\setminus\{o\})$ and for a fixed real number $R>0$, the following distance function
\begin{align*}
	d_C(f,R):= \bigg(\int_{\ba}\frac{|f(x)-\big(\ln \frac{1}{\rho(x)}\big)^{\frac{p-1}{p}}R^{\frac{p-1}{p}}f(e^{-R^{-1}}\sigma)|^p}{\rho^N(x)\big(\ln \frac{1}{\rho(x)}\big)^{p}\big|\ln \big(R\ln \frac{1}{\rho(x)}\big)\big|^{p}}\dv\bigg)^{\frac{1}{p}}.
\end{align*}
Now we have the following stability result in the critical setup:
\begin{theorem}\label{stab-th-crt}
Let $2\leq p <\infty$. Then for all complex-valued functions $f\in C_0^\infty(\ba\setminus\{o\})$, we have 
	\begin{multline*}
		c_p\bigg(\frac{p-1}{p}\bigg)^p \sup_{R>0} d_C(f,R)^p \leq \int_{\ba}\frac{|\partial_\rho f(x)|^p}{\rho^{N-p}(x)}\dv \\- \bigg(\frac{p-1}{p}\bigg)^p\int_{\ba}\frac{|f(x)|^p}{\rho^{N}(x)(\ln \frac{1}{\rho(x)})^p}\dv,
	\end{multline*}
 and applying \eqref{gauss}, we have
\begin{multline*}
		c_p\bigg(\frac{p-1}{p}\bigg)^p \sup_{R>0} d_C(f,R)^p \leq \int_{\ba}\frac{|\nabla_g f(x)|_g^p}{\rho^{N-p}(x)}\dv \\- \bigg(\frac{p-1}{p}\bigg)^p\int_{\ba}\frac{|f(x)|^p}{\rho^{N}(x)(\ln \frac{1}{\rho(x)})^p}\dv,
	\end{multline*}
where $c_p$ is defined in Lemma~\ref{cpl}.   
\end{theorem}

\subsection{Double weighted Hardy inequality}\label{sub-3} Here we present the double-weighted Hardy-type inequality where the singularity occurs in the origin as well as on the boundary. For simplicity, we state the results for the unit ball of the manifold $\ba$, but one can easily mimic the proof on $\br$ for any $R>0$.
\begin{theorem}\label{new-hardy-th}
    Let $1< p<\infty$, $a<N$, $p\leq b<\infty$, $c>0$, and $a\leq N- (b-1)c$. Then for any complex-valued function $f\in C_0^\infty (\ba)$ there holds
\begin{align}\label{new-wg-hardy}
        \bigg(\frac{(b-1)}{p}c\bigg)^p\int_{\ba} \frac{|f(x)|^{p}}{\rho^{a}(x)(1-\rho^{c}(x))^{b}} \dv\leq \int_{\ba} \frac{|\partial_{\rho} f(x)|^p}{\rho^{a-p}(x)(1-\rho^{c}(x))^{b-p}} \dv,
    \end{align}
    and applying \eqref{gauss} we have 
 \begin{align}
        \bigg(\frac{(b-1)}{p}c\bigg)^p\int_{\ba} \frac{|f(x)|^{p}}{\rho^{a}(x)(1-\rho^{c}(x))^{b}} \dv\leq \int_{\ba} \frac{|\nabla_g f(x)|_g^p}{\rho^{a-p}(x)(1-\rho^{c}(x))^{b-p}} \dv,
\end{align}   
where the constant $\big(\frac{(b-1)}{p}c\big)^p$ is sharp.
\end{theorem}

In developing the above result, we obtain new Caffarelli--Kohn--Nirenberg type inequalities in the Cartan--Hadamard manifolds.
\begin{corollary}\label{ckn-th}
Let $1< p, \, q <\infty$, $a<N$, $p\leq b<\infty$, $c>0$, and $a\leq N- (b-1)c$. Let $0<\eta<\infty$ with $p+q\geq \eta$ and $\delta\in [0,1]\cap \big[\frac{\eta-q}{\eta}, \frac{p}{\eta}\big]$ and $\alpha,\beta,\gamma\in \mathbb{R}$. Assume that, $\frac{\delta \eta}{p}+\frac{(1-\delta)\eta}{q}=1$ and $\gamma = \delta (\alpha-1)+\beta (1-\delta)$. Then for any complex-valued function $f\in C_0^\infty (\ba)$ there holds
    \begin{multline}\label{ckn}
        ||w^\gamma(x) f(x)||_{L^\eta(\ba)}\leq \bigg|\frac{p}{c(b-1)}\bigg|^\delta\bigg|\bigg|\frac{\partial_\rho f(x)}{\rho^{\frac{a-p}{p}}(x)(1-\rho^c(x))^{\frac{b-p}{p}}}\bigg|\bigg|^\delta_{L^p(\ba)}\\ \times ||w^\beta(x) f(x)||_{L^q(\ba)}^{1-\delta},
    \end{multline}
    and applying \eqref{gauss}, we have
    \begin{multline}\label{ckn-o}
        ||w^\gamma(x) f(x)||_{L^\eta(\ba)}\leq \bigg|\frac{p}{c(b-1)}\bigg|^\delta\bigg|\bigg|\frac{|\nabla_g f(x)|_g}{\rho^{\frac{a-p}{p}}(x)(1-\rho^c(x))^{\frac{b-p}{p}}}\bigg|\bigg|^\delta_{L^p(\ba)}\\ \times ||w^\beta(x) f(x)||_{L^q(\ba)}^{1-\delta},
    \end{multline}
    where $w(x)=\rho^{\frac{a}{p(1-\alpha)}}(x)(1-\rho^c(x))^{\frac{b}{p(1-\alpha)}}$. The constant in the inequality \eqref{ckn}, and  \eqref{ckn-o} is sharp for $\delta=0$ or $\delta = 1$.
\end{corollary}
\begin{remark}
    For the stability of the Caffarelli--Kohn--Nirenberg (weighted Hardy--Sobolev type) inequality, we can refer to the very recent paper \cite{FP23} and references therein.
\end{remark}

Now by substituting $c=\frac{N-a}{p-1}$, and $b=p$ into Theorem~\ref{new-hardy-th}, we obtain an improvement of known weighted $L^p$-Hardy inequality. This can be compared with \cite[Theorem~3.1]{vhn}. Let us mention this improvement of classical Hardy inequalities in the following remark.
\begin{remark}
 Let $1< p<\infty$, and $a<N$. Then for any complex-valued function $f\in C_0^\infty (\ba)$ there holds
\begin{multline*}
        \int_{\ba} \frac{|f(x)|^{p}}{\rho^{a}(x)} \dv\leq \int_{\ba} \frac{|f(x)|^{p}}{\rho^{a}(x)(1-\rho^{\frac{N-a}{p-1}}(x))^{p}} \dv \\ \leq \bigg(\frac{p}{N-a}\bigg)^p\int_{\ba} \frac{|\partial_\rho f(x)|^p}{\rho^{a-p}(x)} \dv,
    \end{multline*}
and applying \eqref{gauss} we have 
\begin{multline*}
        \int_{\ba} \frac{|f(x)|^{p}}{\rho^{a}(x)} \dv\leq \int_{\ba} \frac{|f(x)|^{p}}{\rho^{a}(x)(1-\rho^{\frac{N-a}{p-1}}(x))^{p}} \dv\\ \leq \bigg(\frac{p}{N-a}\bigg)^p\int_{\ba} \frac{|\nabla_g f(x)|_g^p}{\rho^{a-p}(x)} \dv,
    \end{multline*}   
where the constant $\big(\frac{p}{N-a}\big)^p$ is sharp.
\end{remark}

 We know that $1-t^\alpha=\alpha\ln \frac{1}{t}+ o(\alpha)$ as $\alpha\rightarrow 0$. Using this in Theorem~\ref{new-hardy-th}, we deduce the following result.
\begin{remark}
  Let $1< p<\infty$, $p\leq b<\infty$, and $a\leq N$. Then for any complex-valued function $f\in C_0^\infty (\ba)$ there holds 
    \begin{align*}
        \bigg(\frac{b-1}{p}\bigg)^p\int_{\ba} \frac{|f(x)|^{p}}{\rho^{a}(x)(\ln \frac{1}{\rho(x)})^{b}} \dv\leq \int_{\ba} \frac{|\partial_\rho f(x)|^p}{\rho^{a-p}(x)(\ln \frac{1}{\rho(x)})^{b-p}} \dv,
    \end{align*}
    and applying \eqref{gauss} we have 
    \begin{align*}
        \bigg(\frac{b-1}{p}\bigg)^p\int_{\ba} \frac{|f(x)|^{p}}{\rho^{a}(x)(\ln \frac{1}{\rho(x)})^{b}} \dv\leq \int_{\ba} \frac{|\nabla_g f(x)|_g^p}{\rho^{a-p}(x)(\ln \frac{1}{\rho(x)})^{b-p}} \dv,
    \end{align*}
    where the constant $\big(\frac{b-1}{p}\big)^p$ is sharp. The sharpness can be proved by taking a similar minimizing sequence used in Theorem~\ref{new-hardy-th}. This time we need to multiply cut off function by $\left(\ln\left(\frac{1}{\rho(x)}\right)\right)^\mu$ for some $\mu>\left(\frac{b-1}{p}\right)$.  Moreover, for the cases of $a=N$ and $b=p$, we deduce the known critical Hardy inequality result studied in \cite[Theorem~3.3]{vhn}.
\end{remark}

We conclude this subsection by mentioning an improvement of the geometric Hardy inequality in terms of the distance function $\text{dist}(x,\partial \ba)=1-\rho(x)$. Substituting $a=p,\, c=1$ in Theorem~\ref{new-hardy-th}, we have the following remark.
\begin{remark}
     Let $1< p<\infty$, $p<N$, $p\leq b<\infty$, and $p\leq N-b+1$. Then for any complex-valued function $f\in C_0^\infty (\ba)$ there holds
\begin{align*}
       \int_{\ba} \frac{|f(x)|^{p}}{\text{dist}(x,\partial \ba)^b} \dv &\leq  \int_{\ba} \frac{|f(x)|^{p}}{\rho^{p}(x)(1-\rho(x))^{b}} \dv\\&\leq \bigg(\frac{p}{b-1}\bigg)^p \int_{\ba} \frac{|\partial_\rho f(x)|^p}{\text{dist}(x,\partial \ba)^{b-p}} \dv,
    \end{align*}
    and applying \eqref{gauss} we have 
\begin{align*}
       \int_{\ba} \frac{|f(x)|^{p}}{\text{dist}(x,\partial \ba)^b} \dv &\leq  \int_{\ba} \frac{|f(x)|^{p}}{\rho^{p}(x)(1-\rho(x))^{b}} \dv\\&\leq \bigg(\frac{p}{b-1}\bigg)^p \int_{\ba} \frac{|\nabla_g f(x)|_g^p}{\text{dist}(x,\partial \ba)^{b-p}} \dv.
    \end{align*}
Furthermore, we can notice for any complex-valued function $f\in C_0^\infty (\ba)$ there holds
\begin{align*}
    \int_{\ba} \frac{|f(x)|^{p}}{\text{dist}(x,\partial \ba)^b} \dv \leq  \int_{\ba} \frac{|f(x)|^{p}}{\rho^{p}(x)(1-\rho(x))^{b}} \dv,
\end{align*}
where the constant $1$ is sharp. Therefore, the constant $\big(\frac{p}{b-1}\big)^p$ that appears in the earlier two inequalities is also sharp.
\end{remark}

\subsection{Higher order improvements}\label{sub-4} In this part, we apply the two-weighted Hardy inequality  \eqref{new-wg-hardy} to obtain a new type of the weighted Rellich inequality. Let us recall \cite[Lemma~4.1]{vhn}: Let $(\m,g)$ be an $N$-dimensional Cartan--Hadamard manifold with $N\geq 2$. For any $p\in (1, N)$ and $-N(p-1)<\beta<N-p$, there holds 
\begin{align}\label{rell-typ}
\bigg(\frac{N(p-1)+\beta}{p}\bigg)^p\int_{\m}\frac{|f(x)|^p}{\rho^{p+\beta}(x)}\dv\leq \int_{\m}\frac{\big|\partial_\rho f(x)+\big(\frac{N-1}{\rho(x)}+\frac{J_\rho(\rho(x),\sigma)}{J(\rho(x),\sigma)}\big)f(x)\big|^p}{\rho^\beta(x)}\dv,
\end{align}
for any $f\in C_0^\infty(\m)$. Moreover, the constant $\big(\frac{N(p-1)+\beta}{p}\big)^p$ is the best possible in \eqref{rell-typ}. Combining this result with Theorem~\ref{new-hardy-th}, we have the following two weighted Rellich-type inequalities. In some particular cases, we deduce a different type of improvement of \cite[Theorem~4.3]{vhn}.
\begin{corollary}\label{rel-cor}
Let $1< p<N$, $2p+\beta<N$, $-N(p-1)<\beta<N-p$, $c>0$, and $2p+\beta\leq N- (p-1)c$. Then for any complex-valued function $f\in C_0^\infty (\ba)$ there holds
\begin{multline}\label{rel-cor-eqn-ano}
        \bigg(\frac{N(p-1)+\beta}{p}\bigg)^p\bigg(\frac{(p-1)}{p}c\bigg)^p\int_{\ba} \frac{|f(x)|^{p}}{\rho^{2p+\beta}(x)(1-\rho^{c}(x))^{p}} \dv \\ \leq \int_{\ba} \frac{|\Delta_{\rho,g} f(x)|^p}{\rho^{\beta}(x)} \dv.
\end{multline}
In particular choosing $c=\frac{N-2p-\beta}{p-1}$, for any $f\in C_0^\infty (\ba)$ there holds
\begin{align}\label{rel-cor-eqn}
        \int_{\ba} \frac{|f(x)|^{p}}{\rho^{2p+\beta}(x)} \dv &\leq \int_{\ba} \frac{|f(x)|^{p}}{\rho^{2p+\beta}(x)(1-\rho^{\frac{N-2p-\beta}{p-1}}(x))^{p}} \dv\nonumber\\&\leq \frac{ p^{2p}}{\big(N(p-1)+\beta\big)^p\big(N-2p-\beta\big)^p}\int_{\ba} \frac{|\Delta_{\rho,g} f(x)|^p}{\rho^{\beta}(x)} \dv,
\end{align}
 and both the constants in \eqref{rel-cor-eqn-ano} and \eqref{rel-cor-eqn} are sharp. 
\end{corollary}

Presenting the aforementioned result for the complete Laplacian operator is not feasible at this point, as a comparison between the radial and full Laplacian operators on the Riemannian manifold is not known in general. However, there are instances where such a comparison is available. As an illustration, we can refer to the comparison lemma (as detailed in \cite[Lemma~6.1]{jmaa}). Let $(\m,g,\psi)$, be an $N$-dimensional Riemannian model manifold with $\psi$ satisfying the sectional curvature condition $-\frac{\psi^{\prime\prime}(\rho(x))}{\psi(\rho(x))}\leq -1$. Then, for $0\leq \beta <N-4$, there holds
\begin{align*}
    \int_{\m}\frac{|\Delta_{\rho,g} f(x)|^2}{\rho^\beta(x)}\dv\leq \int_{\m}\frac{|\Delta_{g} f(x)|^2}{\rho^\beta(x)}\dv,
\end{align*}
for all $f\in C_0^\infty(\m)$. Moreover, equality holds when $f$ is a radial function. Looking at the proof of this observation closely, one can notice that the above comparison is independent of the domain. Hence, taking $p=2$ in the Corollary~\ref{rel-cor} and with the help of comparison inequality, we have the following remark.
\begin{remark}\label{rem-38}
  Let $(\m,g,\psi)$ be an $N$-dimensional Riemannian model manifold with $N\geq 4$ and $-\frac{\psi^{\prime\prime}(\rho(x))}{\psi(\rho(x))}\leq -1$. Let $\ba\subset\m$ be the unit ball. Let $c>0$, and $0\leq \beta \leq N-4-c$. Then for any complex-valued function $f\in C_0^\infty (\ba)$ there holds
\begin{align*}
        \bigg(\frac{N+\beta}{4}c\bigg)^2\int_{\ba} \frac{|f(x)|^{2}}{\rho^{4+\beta}(x)(1-\rho^{c}(x))^{2}} \dv\leq \int_{\ba} \frac{|\Delta_{g} f(x)|^2}{\rho^{\beta}(x)} \dv.
\end{align*}
In particular, let us choose  $c=(N-\beta-4)$, and using $1\leq \big(1-\rho^{N-\beta-4}(x)\big)^{-1}$ on $\ba$, we deduce
\begin{align*}
        \int_{\ba} \frac{|f(x)|^{2}}{\rho^{4+\beta}(x)} \dv & \leq  \int_{\ba} \frac{|f(x)|^{2}}{\rho^{4+\beta}(x)(1-\rho^{N-\beta-4}(x))^{2}} \dv \\&\leq \frac{16}{(N+\beta)^2(N-\beta-4)^2} \int_{\ba} \frac{|\Delta_{g} f(x)|^2}{\rho^{\beta}(x)} \dv.
\end{align*}
Here, above the constant is sharp and we get a different type of improved version of the known Rellich inequality on $\ba$. For the flat manifold (i.e., $\psi(\rho(x))=\rho(x)$) case, it is worth mentioning that the required curvature condition does not hold, and to get this type of improvement, refer to \cite[Remark~1.14]{SY-arxiv}.
\end{remark}

We end this section by mentioning the higher-order version of this result for the higher-order operator introduced in the preliminaries section. By choosing an appropriate $c$, this result can be compared with \cite[Theorem~4.5]{vhn}. Also using the fact $1-t^\alpha=\alpha\ln \frac{1}{t}+ o(\alpha)$ as $\alpha\rightarrow 0^+$, one can compare this result with \cite[Theorem~4.6]{vhn}. The work in \cite[Example~2, p.~879]{Barbatis} presents improvements to $L^p$  Rellich inequalities of any order, incorporating sharp logarithmic terms, for arbitrary Cartan--Hadamard manifolds. In the Euclidean case, we refer to \cite[Theorem~2.3]{adi} for Hardy--Rellich inequalities involving polyharmonic operators in the critical dimension, which include logarithmic terms. 

We define for $\beta\, , N$, and $p$ the following notation
\begin{align*}
    \Lambda_{i}(N,p,\beta)= \bigg(\frac{N(p-1)+\beta+ip}{p}\bigg)^p \text{ for integer } i\geq 0 \quad \text{ and } \Lambda_{-1}(N,p,\beta)=1.
\end{align*}
Now we are ready to present our result.
\begin{theorem}\label{higher-order}
Let $(\m,g)$ be an $N$-dimensional Cartan–-Hadamard manifold with $N\geq 2$. Let $\ba\subset\m$ be the unit ball and $k\geq 1$ be an integer. Let $1< p<N$, $kp+\beta<N$, $-N(p-1)<\beta<N-(k-1)p$, $c>0$, and $kp+\beta\leq N- (p-1)c$. Let $f\in C_0^\infty (\ba)$ be any complex-valued function. Then for some integer $\alpha$, if $k=2\alpha-1$, there holds
\begin{multline*}
        \bigg(\frac{(p-1)}{p}c\bigg)^{\alpha p} \prod_{i=0}^{\alpha-1} \Lambda_{2i-1}(N,p,\beta)\int_{\ba} \frac{|f(x)|^{p}}{\rho^{kp+\beta}(x)(1-\rho^{c}(x))^{p}} \dv\\ \leq \int_{\ba} \frac{|\partial_\rho \Delta_{\rho,g}^{\alpha-1} f(x)|^p}{\rho^{\beta}(x)} \dv,
\end{multline*}
and if $k=2\alpha$, there holds
\begin{multline*}
        \bigg(\frac{(p-1)}{p}c\bigg)^{\alpha p} \prod_{i=0}^{\alpha-1} \Lambda_{2i}(N,p,\beta) \int_{\ba} \frac{|f(x)|^{p}}{\rho^{kp+\beta}(x)(1-\rho^{c}(x))^{p}} \dv \\ \leq \int_{\ba} \frac{|\Delta_{\rho,g}^\alpha f(x)|^p}{\rho^{\beta}(x)} \dv.
\end{multline*}  
\end{theorem}

\medspace

\section{Proofs of Theorem~\ref{anoth-crit-th} and Corollary~\ref{anoth-crit-cor-rem}}\label{pfs-1}
We will demonstrate the proof of Theorem~\ref{anoth-crit-th} and Corollary~\ref{anoth-crit-cor-rem} in this section.

{\bf Proof of Theorem~\ref{anoth-crit-th}:} Let us start the computation on the geodesic ball $\br$ for some $R>0$ and for notational economy assume $r=\rho(x)$. Using the polar coordinate, we deduce the following
    \begin{align*}
        \int_{\br}&\frac{|f-f_R|^p}{|r^{\frac{N}{p}}\big(\ln\frac{R}{r}\big)^{\frac{\gamma}{p}}|^p}\dv=\int_{0}^R\int_{\sn}\frac{|f(r\sigma)-f(R\sigma)|^p}{r^{N}\big(\ln\frac{R}{r}\big)^{\gamma}}r^{N-1}J(r,\sigma)\dsn\dr\\&=\int_{0}^R\int_{\sn}\frac{|f(r\sigma)-f(R\sigma)|^p}{\big(\ln\frac{R}{r}\big)^{\gamma}}r^{-1}J(r,\sigma)\dsn\dr\\&=\int_{0}^R \frac{d}{dr}\bigg(\int_{\sn}\frac{|f(r\sigma)-f(R\sigma)|^p}{(\gamma-1)\big(\ln\frac{R}{r}\big)^{\gamma-1}}J(r,\sigma)\dsn\bigg)\dr\\&-\frac{p}{\gamma-1}{\rm Re}\int_{0}^R \int_{\sn}\frac{|f(r\sigma)-f(R\sigma)|^{p-2}(f(r\sigma)-f(R\sigma))}{\big(\ln\frac{R}{r}\big)^{\gamma-1}}\overline{\frac{df(r\sigma)}{dr}}J(r,\sigma)\dsn\dr\\&-\frac{1}{\gamma-1}\int_{0}^R \int_{\sn}\frac{|f(r\sigma)-f(R\sigma)|^p}{\big(\ln\frac{R}{r}\big)^{\gamma-1}}J_r(r,\sigma)\dsn\dr.
    \end{align*}
    Now notice that $|f(r\sigma)-f(R\sigma)|\leq \|\grad_g f\|_{L^\infty(M)}(R-r)\leq R\|\grad_g f\|_{L^\infty(M)} \ln (R/r)$ for $0<r\leq R$. Then, using this along with $p>\gamma-1$, the boundary term in the first integral vanishes at $r=R$. For the other boundary in the first integral condition at $r=0$, use $|f(r\sigma)-f(R\sigma)|\leq \|\grad_g f\|_{L^\infty(M)}(R-r)$ and $\gamma>1$. This gives the first integral vanishes. The last term is always nonpositive due to the fact \eqref{den-der} on Cartan--Hadamard manifolds, and so $J_r(r,\sigma)\geq 0$. So we deduce
    \begin{align}\label{steps-ref}
        \nonumber&\int_{\br}\frac{|f-f_R|^p}{|r^{\frac{N}{p}}\big(\ln\frac{R}{r}\big)^{\frac{\gamma}{p}}|^p}\dv=\int_{0}^R\int_{\sn}\frac{|f(r\sigma)-f(R\sigma)|^p}{r^{N}\big(\ln\frac{R}{r}\big)^{\gamma}}r^{N-1}J(r,\sigma)\dsn\dr\nonumber\\&\leq -\frac{p}{(\gamma-1)}{\rm Re}\int_{0}^R \int_{\sn}\frac{|f(r\sigma)-f(R\sigma)|^{p-2}(f(r\sigma)-f(R\sigma))}{\big(\ln\frac{R}{r}\big)^{\gamma-1}}\overline{\frac{df(r\sigma)}{dr}}J(r,\sigma)\dsn\dr\nonumber\\&\leq \frac{p}{(\gamma-1)}\int_{0}^R \int_{\sn}\frac{|f(r\sigma)-f(R\sigma)|^{p-1}}{\big(\ln\frac{R}{r}\big)^{\gamma-1}}\bigg|\frac{df(r\sigma)}{dr}\bigg|J(r,\sigma)\dsn\dr \nonumber\\&=\frac{p}{(\gamma-1)}\int_{\br}\frac{|f-f_R|^{p-1}}{r^{N-1}\big(\ln\frac{R}{r}\big)^{\gamma-1}}\bigg|\frac{df}{dr}\bigg|\dv \nonumber\\&\leq \frac{p}{(\gamma-1)}\bigg(\int_{\br}\frac{|f-f_R|^{p}}{\big|r^{\frac{N}{p}}\big(\ln\frac{R}{r}\big)^{\frac{\gamma}{p}}\big|^p}\dv\bigg)^{(p-1)/p} \bigg(\int_{\br}\bigg|\frac{df}{dr}\bigg|^pr^{(p-N)}\bigg(\ln\frac{R}{r}\bigg)^{(p-\gamma)}\dv\bigg)^{1/p}.
    \end{align}
    Thus, we obtain 
    \begin{align*}
        \bigg(\int_{\br}\frac{|f-f_R|^p}{|r^{\frac{N}{p}}\big(\ln\frac{R}{r}\big)^{\frac{\gamma}{p}}|^p}\dv\bigg)^{\frac{1}{p}}\leq \frac{p}{(\gamma-1)}\bigg(\int_{\br}r^{(p-N)}\bigg|\ln\frac{R}{r}\bigg|^{(p-\gamma)}\bigg|\frac{df}{dr}\bigg|^p\dv\bigg)^{\frac{1}{p}}.
    \end{align*}
       Similarly, one can observe that
    \begin{align*}
        \bigg(\int_{\br^c}\frac{|f-f_R|^p}{|r^{\frac{N}{p}}\big(\ln\frac{R}{r}\big)^{\frac{\gamma}{p}}|^p}\dv\bigg)^{\frac{1}{p}}\leq \frac{p}{(\gamma-1)}\bigg(\int_{\br^c}r^{(p-N)}\bigg|\ln\frac{R}{r}\bigg|^{(p-\gamma)}\bigg|\frac{df}{dr}\bigg|^p\dv\bigg)^{\frac{1}{p}}.
    \end{align*}
    Combining both, we deduce the required inequality. 
    
  Now we verify the optimality of the constant $\frac{(\gamma-1)}{p}$. First note that from \eqref{anoth-crit-th}
  it implies     \begin{align}\label{help}
        \bigg(\int_{\br}\frac{|f(r\sigma)|^p}{r^{N}\big(\ln\frac{R}{r}\big)^{\gamma}}\dv\bigg)^{\frac{1}{p}}\leq \frac{p}{(\gamma-1)}\bigg(\int_{\br}r^{(p-N)}\bigg(\ln\frac{R}{r}\bigg)^{(p-\gamma)}\bigg|\frac{df}{dr}\bigg|^p\dv\bigg)^{\frac{1}{p}}
    \end{align}
    for all $f\in C_0^\infty(\br)$.
    Therefore, it is enough to consider optimality for the above inequality. Suppose $\delta$ is some small positive number less than $R$, and then consider the sequence of functions
\begin{equation*}
 f_\delta(x)=f_\delta(r\sigma):=
\begin{dcases}
\bigg(\ln\frac{R}{\delta}\bigg)^{\frac{\gamma-1}{p}} &  \text{ when } r\leq \delta; \\
\bigg(\ln\frac{R}{r}\bigg)^{\frac{\gamma-1}{p}} & \text{ when } \delta\leq r\leq  \frac{R}{2}; \\
\left(\ln 2 \right)^{\frac{\gamma-1}{p}}\frac{2}{R}(R-r) &  \text{ when }\frac{R}{2}\leq r\leq  R. \\
\end{dcases}
\end{equation*}
This gives that 
\begin{equation*}
\frac{df_\delta}{dr}(r\sigma)=
\begin{dcases}
0 &  \text{ when } r<\delta; \\
-\left(\frac{\gamma-1}{p}\right)\bigg(\ln\frac{R}{r}\bigg)^{\frac{\gamma-1}{p}-1}r^{-1} & \text{ when } \delta < r <  \frac{R}{2}; \\
-\left(\ln 2 \right)^{\frac{\gamma-1}{p}}\frac{2}{R} & \text{ when }\frac{R}{2} < r <  R. \\
\end{dcases}
\end{equation*}
We can check $f_\delta\in W^{1,p}_{0}(\br)$ and it is enough to consider these as test functions because \eqref{help} holds on this Sobolev space due to the density arguments. Now, from polar coordinate decomposition and increasing nature of $J(r,\sigma)$, and from \eqref{imp-den} near small $\delta$ it follows that
    \begin{align*}
       \int_{\br}\frac{|f_\delta(r\sigma)|^p}{r^{N}\big(\ln\frac{R}{r}\big)^{\gamma}}\dv &\geq\int_{\mathcal{B}_{\frac{R}{2}}(o)\setminus \mathcal{B}_\delta(o)}\frac{|f_\delta(r\sigma)|^p}{r^{N}\big(\ln\frac{R}{r}\big)^{\gamma}}\dv\\&=  \int_{\sn}\int_\delta^{\frac{R}{2}} \bigg(\ln\frac{R}{r}\bigg)^{-1}\frac{J(r,\sigma)}{r} \dr\dsn \\& \geq |\sn|\int_\delta^\frac{R}{2} \bigg(\ln\frac{R}{r}\bigg)^{-1}\frac{1}{r} \dr\\&=|\sn|\left(\ln\left(\ln\frac{R}{\delta}\right)-\ln(\ln2))\right).
    \end{align*}
    So, we have 
    \begin{align*}
        \lim_{\delta\rightarrow 0^+} \int_{\br}\frac{|f_\delta(r\sigma)|^p}{r^{N}\big(\ln\frac{R}{r}\big)^{\gamma}}\dv = +\infty.
    \end{align*}    
    Also, we deduce
    \begin{align*}
        &\int_{\br}r^{(p-N)}\left(\ln\frac{R}{r}\right)^{(p-\gamma)}\left|\frac{df_\delta}{dr}\right|^p\dv\\&\leq \bigg(\frac{\gamma-1}{p}\bigg)^p\int_{\mathcal{B}_{\frac{R}{2}}(o)\setminus \mathcal{B}_\delta(o)}\frac{|f_\delta(r\sigma)|^p}{r^{N}\big(\ln\frac{R}{r}\big)^{\gamma}}\dv\\&+(\ln 2)^{\gamma-1}2^pR^{-p}\int_{\mathcal{B}_R(o)\setminus \mathcal{B}_{\frac{R}{2}}(o)}r^{(p-N)}\left(\ln\frac{R}{r}\right)^{(p-\gamma)}\dv.
    \end{align*}
    We also deduce there is some positive constant $c$, which depends on the bound of density function $J(r,\sigma)$ on $\mathcal{B}_R(o)\setminus \mathcal{B}_{\frac{R}{2}}(o)$ and there holds
    \begin{align*}
        \int_{\mathcal{B}_R(o)\setminus \mathcal{B}_{\frac{R}{2}}(o)}r^{(p-N)}\left(\ln\frac{R}{r}\right)^{(p-\gamma)}\dv&\leq cR^{p}|\sn|\int_{\frac{R}{2}}^R\left(\ln\frac{R}{r}\right)^{(p-\gamma)}\frac{1}{r}\dr\\&=\frac{cR^{p}|\sn|}{(p-\gamma+1)}\left(\ln2\right)^{(p-\gamma+1)}<+\infty,
    \end{align*}
    where we used the fact that $p-\gamma+1>0$. This gives in combining that
    \begin{align*}
       \lim_{\delta\rightarrow 0^+}\frac{\int_{\br}r^{(p-N)}\big(\ln\frac{R}{r}\big)^{(p-\gamma)}\big|\frac{df_\delta}{dr}\big|^p\dv}{\int_{\br}\frac{|f_\delta(r\sigma)|^p}{r^{N}\big(\ln\frac{R}{r}\big)^{\gamma}}\dv} =\bigg(\frac{\gamma-1}{p}\bigg)^p.
    \end{align*}

{\bf Proof of Corollary~\ref{anoth-crit-cor-rem}:} This follows from the steps of Theorem~\ref{anoth-crit-th}. Applying the H\"older inequality after combining the two integrals, one on $\br$ and another on $\br^c$, we deduce an improvement of \eqref{anoth-crit} and \eqref{ful-anoth-crit} with remainder terms as desired.

\medspace

\section{Proofs of Theorem~\ref{new_critical_hardy} and Corollary~\ref{new-critical-hardy-rem}}\label{pfs-2}
This section is devoted to the proofs of Theorem~\ref{new_critical_hardy} and Corollary~\ref{new-critical-hardy-rem}. 

{\bf Proof of Theorem~\ref{new_critical_hardy}:} Let supp$(f)\subset \br\setminus\{o\}$ and assume $r=\rho(x)$. Then, using polar coordinate decomposition and exploiting integration by parts, we obtain
\begin{equation}\label{pr_thm3.2_1}
\begin{split}
  & \int_{\m}\frac{|f|^p}{r^N}\dv\\& =  \int_{\sn}\int_{0}^{R}|f(\text{Exp}_o(r\sigma))|^p\:J(r,\sigma)r^{-1}\dr\dsn\\&= \int_{\sn}\int_{0}^{R}|F(r\sigma)|^p\:J(r,\sigma)\,{\rm d}(\ln r)\dsn\\&= -p{\rm Re}\int_{\sn}\int_{0}^{R}|F|^{p-2}F \overline{F_r}\:J(r,\sigma)(\ln r)\dr \dsn -\int_{\sn}\int_{0}^{R}|F|^{p}\:J_r(r,\sigma)(\ln r)\dr \dsn\\&= -p{\rm Re}\int_{\m}(\ln r) r^{1-N} |f|^{p-2}\, f \, \overline{\frac{\partial f}{\partial r}} \, \dv- \int_{\m} |f|^p \frac{J_r(r,\sigma)}{J(r,\sigma)}\, (\ln r) r^{1-N}\dv.
\end{split}
\end{equation}
In the above, we used
\begin{align*}
    \left[\int_{\sn}|F(r\sigma)|^pJ(r,\sigma)(\ln r)\dsn\right]_{r=0}^{r=R}=0,
\end{align*}
because of \eqref{imp-den} and $f$ is smooth and compactly supported inside $\br\setminus\{o\}$. Now using the representation \eqref{Rp_1} we have
\begin{multline*}
    \int_{\m} R_p\bigg(r^{-\frac{N}{p}}f,-p\frac{(\ln r)}{r^{\frac{(N-p)}{p}}}\frac{\partial f}{\partial r}\bigg)\left|r^{-\frac{N}{p}}f+p\frac{(\ln r)}{r^{\frac{(N-p)}{p}}}\frac{\partial f}{\partial r}\right|^{2}\dv \\ = \frac{p^p}{p} \int_{\m} \frac{|\ln r|^p}{r^{N-p}}\bigg|\frac{\partial f}{\partial r}\bigg|^p\dv + \frac{(p-1)}{p}\int_{\m}\frac{|f|^p}{r^{N}}\dv +p{\rm Re}\int_{\m} |f|^{p-2}\, f \,\overline{\frac{\partial f}{\partial r}} \, (\ln r)\, r^{1-N}\dv.
\end{multline*}
Therefore, we deduce
\begin{multline*}
    \int_{\m}\frac{|f|^p}{r^N}\dv  = -p{\rm Re}\int_{\m} |f|^{p-2}\, f \, \overline{\frac{\partial f}{\partial r}} \, (\ln r) r^{1-N}\dv- \int_{\m} |f|^p \frac{J_r(r,\sigma)}{J(r,\sigma)}\, (\ln r) r^{1-N}\dv \\ = -\int_{\m} R_p\bigg(r^{-\frac{N}{p}}f,-p\frac{(\ln r)}{r^{\frac{(N-p)}{p}}}\frac{\partial f}{\partial r}\bigg)\left|r^{-\frac{N}{p}}f+p\frac{(\ln r)}{r^{\frac{(N-p)}{p}}}\frac{\partial f}{\partial r}\right|^{2}\dv \\+\frac{p^p}{p} \int_{\m} \frac{|\ln r|^p}{r^{N-p}}\bigg|\frac{\partial f}{\partial r}\bigg|^p\dv+\frac{(p-1)}{p}\int_{\m}\frac{|f|^p}{r^N}\dv- \int_{\m} |f|^p \frac{J_r(r,\sigma)}{J(r,\sigma)}\, (\ln r) r^{1-N}\dv.
\end{multline*}
This implies that
\begin{multline}\label{last_ineq}
   \int_{\m}\frac{|f|^p}{r^N}\dv=p^p \int_{\m} \frac{|\ln r|^p}{r^{N-p}}\bigg|\frac{\partial f}{\partial r}\bigg|^p\dv - p\int_{\m} |f|^p \frac{J_r(r,\sigma)}{J(r,\sigma)}\, (\ln r) r^{1-N}\dv \\-p\int_{\m} R_p\bigg(r^{-\frac{N}{p}}f,-p\frac{(\ln r)}{r^{\frac{(N-p)}{p}}}\frac{\partial f}{\partial r}\bigg)\left|r^{-\frac{N}{p}}f+p\frac{(\ln r)}{r^{\frac{(N-p)}{p}}}\frac{\partial f}{\partial r}\right|^{2}\dv,
\end{multline}
which is \eqref{new_critical_hardy_1}.

Now, to show \eqref{new_critical_hardy_2} we use the H\"older inequality in \eqref{pr_thm3.2_1} to get
\begin{align*}
  & \int_{\m}\frac{|f|^p}{r^N}\dv\\& =  -p{\rm Re}\int_{\m}(\ln r) r^{1-N} |f|^{p-2}\, f \, \overline{\frac{\partial f}{\partial r}} \, \dv- \int_{\m} |f|^p \frac{J_r(r,\sigma)}{J(r,\sigma)}\, (\ln r) r^{1-N}\dv \\&
  \leq p\left(\int_{\m}\left|\frac{|f|^{p-2}f}{r^{\frac{N(p-1)}{p}}}\right|^{\frac{p}{p-1}}\dv\right)^{\frac{p-1}{p}}\left(\int_{\m}\left|\frac{\partial f}{\partial r}\right|^{p}\frac{|\ln r|^{p}}{r^{N-p}}\dv\right)^{\frac{1}{p}}\\& \quad \quad - \int_{\m} |f|^p \frac{J_r(r,\sigma)}{J(r,\sigma)}\, (\ln r) r^{1-N}\dv,
\end{align*}
that is,
\begin{multline}
  \left(\int_{\m}\frac{|f|^p}{r^N}\dv\right)^{\frac{1}{p}} \leq p\left(\int_{\m}\left|\frac{\partial f}{\partial r}\right|^{p}\frac{|\ln r|^{p}}{r^{N-p}}\dv\right)^{\frac{1}{p}}\\- \left(\int_{\m} |f|^p \frac{J_r(r,\sigma)}{J(r,\sigma)}\, (\ln r) r^{1-N}\dv\right) \left(\int_{\m}\frac{|f|^p}{r^N}\dv\right)^{\frac{1-p}{p}}.
\end{multline}
{\bf Proof of Corollary~\ref{new-critical-hardy-rem}:} Take $r = \rho(x)$. Since $\m$ satisfies \eqref{condi}, we have $J_r(r, \sigma) \equiv 0$. Using this fact, we obtain \eqref{complex_p2_thm3.2_3} and \eqref{complex_p2_thm3.2_4} from \eqref{new_critical_hardy_2} and \eqref{new_critical_hardy_3}, respectively.

We now prove the optimality of the constant $p^p$ by constructing a minimizing sequence.  For each $\delta$ with $0<\delta<\min\{\frac{1}{4},\frac{1}{p}\}$, we define the radially symmetric function $\{f_\delta\}$ as follows
\begin{equation*}
 f_\delta(r):=
\begin{dcases}
0 &  \text{ when } 0\leq r\leq  \delta; \\
4\left(\ln 4 \right)^{-\frac{1}{p}+\delta}\frac{(r-\delta)}{(1-4\delta)} & \text{ when } \delta\leq r\leq  \frac{1}{4}; \\
\left(-\ln r\right)^{-\frac{1}{p}+\delta} & \text{ when }  \frac{1}{4} \leq r\leq  \frac{1}{1+\delta};\\
\frac{(1+\delta)}{\delta}\left(\ln (1+\delta) \right)^{-\frac{1}{p}+\delta}(1-r) &  \text{ when }\frac{1}{1+\delta}\leq r\leq  1; \\
0 &  \text{ when } 1\leq r <  \infty.
\end{dcases}
\end{equation*}
and so we have
\begin{equation*}
 \frac{\partial f_\delta}{\partial r}=
\begin{dcases}
0 &  \text{ when } 0< r<  \delta; \\
\frac{4}{(1-4\delta)}\left(\ln 4 \right)^{-\frac{1}{p}+\delta} & \text{ when } \delta< r<  \frac{1}{4}; \\
-\left(-\frac{1}{p}+\delta\right)\left(-\ln r\right)^{-\frac{1}{p}+\delta-1}\frac{1}{r} & \text{ when }  \frac{1}{4} < r<  \frac{1}{1+\delta};\\
-\frac{(1+\delta)}{\delta}\left(\ln (1+\delta) \right)^{-\frac{1}{p}+\delta} &  \text{ when }\frac{1}{1+\delta}< r<  1; \\
0 &  \text{ when } 1< r <  \infty.
\end{dcases} 
\end{equation*}
Now, because of $f_\delta(0)=0$ and using boundedness of the function and its derivative on compact support, we can verify that for each $\delta$ we have $f_\delta\in W_0^{1,p}(\m \setminus \{0\})$. Now using (increasing) monotonicity of $J(r,\sigma)$, we have
\begin{align*}
    \int_{\m} \frac{|f_\delta|^p}{r^N}\dv &\geq \left(\int_{\sn}J(1/4,\sigma)\dsn\right)\int_{1/4}^{1/(1+\delta)}(-\ln r)^{-1+\delta p}r^{-1}\dr\\&=\left(\int_{\sn}J(1/4,\sigma)\dsn\right)\left[\frac{(\ln 4)^{\delta p}-(\ln (1+\delta))^{\delta p}}{\delta p}\right].
\end{align*}
Then, using the fact 
\begin{align*}
    \lim_{\delta \rightarrow 0^+}\frac{(\ln 4)^{\delta p}-(\ln (1+\delta))^{\delta p}}{\delta p}=+\infty, 
\end{align*}
we deduce
\begin{align*}
    \int_{\m} \frac{|f_\delta|^p}{r^N}\dv \rightarrow +\infty \text{ as } \delta\rightarrow 0^+.
\end{align*}
On the other hand, it follows that
\begin{align*}
    \int_{\m} \frac{|\ln r|^p\left|\frac{\partial f_\delta}{\partial r}\right|^p}{r^{N-p}}\dv & = \left(\frac{1}{p}-\delta\right)^{p} \int_{ \ball{\frac{1}{1+\delta}}\setminus \ball{\frac{1}{4}}}\frac{|\ln r|^{-1+\delta p}}{r^{N}}\dv\\&+\frac{4^p}{(1-4\delta)^p}\left(\ln 4 \right)^{-1+\delta p} \int_{\ball{\frac{1}{4}}\setminus\ball{\delta}}\frac{|\ln r|^{p}}{r^{N-p}}\dv\\&+\frac{(1+\delta)^p}{\delta^p}\left(\ln (1+\delta) \right)^{-1+\delta p} \int_{\ball{1}\setminus\ball{\frac{1}{1+\delta}}}\frac{|\ln r|^{p}}{r^{N-p}}\dv\\&=\left(\frac{1}{p}-\delta\right)^{p} \int_{ \ball{\frac{1}{1+\delta}}\setminus \ball{\frac{1}{4}}}\frac{|f_\delta|^{p}}{r^{N}}\dv+F_1+F_2,
\end{align*}
where
\begin{align*}
  F_1:=\frac{4^p}{(1-4\delta)^p}\left(\ln 4 \right)^{-1+\delta p} \int_{\ball{\frac{1}{4}}\setminus \ball{\delta}}\frac{|\ln r|^{p}}{r^{N-p}}\dv,
\end{align*}
and
\begin{align*}
 F_2:=\frac{(1+\delta)^p}{\delta^p}\left(\ln (1+\delta) \right)^{-1+\delta p} \int_{\ball{1}\setminus \ball{\frac{1}{1+\delta}}}\frac{|\ln r|^{p}}{r^{N-p}}\dv.
\end{align*}

Next using the fact $p>1$ and $\lim_{\delta \rightarrow 0^+}\delta^{p-1}(-\ln \delta)^p=0$, we deduce
\begin{align*}
    F_1\leq  \frac{4^p}{(1-4\delta)^p}\left(\ln 4 \right)^{-1+\delta p} \left(\int_{\sn}J(1/4,\sigma)\dsn\right)\int_{\delta}^{1/4}r^{p-1}(-\ln r)^p\dr = O(1),
\end{align*}
for $\delta \rightarrow 0^+$. Also using $r\leq 1$ and observing $ \lim_{\delta \rightarrow 0^+}\frac{(1+\delta)^p}{\delta^p}\left(\ln (1+\delta) \right)^{p+\delta p}=1$, we deduce
\begin{align*}
    F_2\leq  \frac{(1+\delta)^p}{\delta^p}\left(\ln (1+\delta) \right)^{-1+\delta p} \left(\int_{\sn}J(1,\sigma)\dsn\right)\int_{\frac{1}{1+\delta}}^{1}r^{-1}(-\ln r)^p\dr=O(1),
\end{align*}
for $\delta \rightarrow 0^+$. Hence, we conclude that 
\begin{align*}
\lim_{\delta \rightarrow 0^+}
\frac{\int_{\m} \frac{|\ln \rho(x)|^p}{\rho^{N-p}(x)}\big|\partial_\rho f_\delta(x)\big|^p\dv}
{\int_{\m} \frac{|f_\delta(x)|^p}{\rho^N(x)}\dv}
= \frac{1}{p^p}.
\end{align*}

\medspace

\section{Proofs of Theorem~\ref{stab-th} and Theorem~\ref{stab-th-crt}}\label{pfs-3}
This section is mainly devoted to proving the stability results. First, let us mention the following result (see \cite[Formula~(2.13)]{FS08}).
\begin{lemma}\label{cpl}
Let $p \geq 2$. Then there exists $c_p > 0$ such that
\begin{align*}
    |a-b|^p \geq |a|^p - p|a|^{p-2}{\rm Re\;} \overline{a}\cdot b + c_p|b|^p
\end{align*}
holds for all vectors $a, b \in \mathbb{C}^N$, where 
\begin{align*}
	c_p=\min_{0<t< 1/2}[(1-t)^p-t^p+pt^{p-1}]
\end{align*}
is sharp in this inequality.
\end{lemma}

{\bf Proof of Theorem~\ref{stab-th}:} Let $x=r\sigma$, where the radial part is $r=\rho(x)$ and the spherical part is described by $\sigma=x/|x|$. For any complex-valued function $f\in C_0^\infty(\m\setminus\{o\})$, we have $f(o)=
\underset{r\rightarrow 0}{\rm lim} f(r\sigma)=0$ and $\underset{r\rightarrow \infty}{\rm lim} f(r\sigma)=0$. Next define the function $g(x)=g(r\sigma):=r^{\frac{N-p-\beta}{p}}f(r\sigma)=|x|^{\frac{N-p-\beta}{p}}f(x)$ with the same boundary conditions as $f$. Also this gives $g\in C_0^\infty(\m\setminus\{o\})$. By using the polar coordinate decomposition, we obtain 
    \begin{equation}
    \begin{split}
    J(f):=&\int_{\m}\frac{|\nabla_{r,g} f|^p}{r^\beta}\dv - \bigg(\frac{N-p-\beta}{p}\bigg)^p\int_{\m}\frac{|f|^p}{r^{p+\beta}}\dv\\&=\int_{\sn}\int_0^\infty\bigg[\bigg|\frac{\partial f}{\partial r}\bigg|^p r^{N-1-\beta}-\bigg(\frac{N-p-\beta}{p}\bigg)^p|f|^p r^{N-1-p-\beta}\bigg]J(r,\sigma)\dr\dsn\\&=\int_{\sn}\int_0^\infty\bigg|\frac{\partial g}{\partial r} r^{-\frac{N-p-\beta}{p}}-\bigg(\frac{N-p-\beta}{p}\bigg)g r^{-\frac{N-\beta}{p}}\bigg|^pr^{N-1-\beta}J(r,\sigma)\dr\dsn\\&-\bigg(\frac{N-p-\beta}{p}\bigg)^p\int_{\sn}\int_0^\infty|f|^p r^{N-1-p-\beta}J(r,\sigma)\dr\dsn.
    \end{split}
\end{equation}

Now applying Lemma~\ref{cpl}, we have 
\begin{align*}
    &\bigg|\frac{\partial g}{\partial r} r^{-\frac{N-p-\beta}{p}}-\bigg(\frac{N-p-\beta}{p}\bigg)g r^{-\frac{N-\beta}{p}}\bigg|^p\\&=\bigg|\bigg(\frac{N-p-\beta}{p}\bigg)g r^{-\frac{N-\beta}{p}}-\frac{\partial g}{\partial r} r^{-\frac{N-p-\beta}{p}}\bigg|^p\\&\geq \bigg(\frac{N-p-\beta}{p}\bigg)^p|g|^pr^{-N+\beta}-p\bigg(\frac{N-p-\beta}{p}\bigg)^{p-1}|g|^{p-2}{\rm Re\;}\overline{g}\frac{\partial g}{\partial r} r^{-N+\beta+1}\\&+c_p\bigg|\frac{\partial g}{\partial r}\bigg|^pr^{-N+p+\beta}\\&= \bigg(\frac{N-p-\beta}{p}\bigg)^p|f|^pr^{-p}-p\bigg(\frac{N-p-\beta}{p}\bigg)^{p-1}|g|^{p-2}{\rm Re\;}\overline{g}\frac{\partial g}{\partial r} r^{-N+\beta+1}+c_p\bigg|\frac{\partial g}{\partial r}\bigg|^pr^{-N+p+\beta}.
\end{align*}
Integration by parts and boundary conditions of the compactly supported smooth function $g$ gives 
\begin{align*}
&\int_{\sn}\int_0^\infty|g|^{p-2}{\rm Re\;}\overline{g}\frac{\partial g}{\partial r} r^{-N+\beta+1}r^{N-1-\beta} J(r,\sigma)\dr\dsn\\&=\int_{\sn}\int_0^\infty\frac{\partial |g|^p}{\partial r} J(r,\sigma)\dr\dsn\\&=-\frac{1}{p}\int_{\sn}\int_{0}^\infty |g|^pJ_r(r,\sigma)\dr\dsn .
\end{align*}
Therefore, we deduce by using \eqref{den-der}, the remainder term as below
\begin{align*}
    J(f)&\geq c_p\int_{\sn}\int_0^
    \infty \bigg|\frac{\partial g}{\partial r}\bigg|^pr^{-N+p+\beta}r^{N-1-\beta}J(r,\sigma)\dr\dsn\\&+\bigg(\frac{N-p-\beta}{p}\bigg)^{p-1}\int_{\sn}\int_{0}^\infty |g|^pJ_r(r,\sigma)\dr\dsn \\&=  c_p\int_{\sn}\int_0^
    \infty \bigg|\frac{\partial g}{\partial r}\bigg|^pr^{p-N}r^{N-1}J(r,\sigma)\dr\dsn\\&+\bigg(\frac{N-p-\beta}{p}\bigg)^{p-1}\int_{\sn}\int_{0}^\infty |g|^pr^{1-N}\frac{J_r(r,\sigma)}{J(r,\sigma)}r^{N-1}J(r,\sigma)\dr\dsn \\& \geq c_p\int_{\m}\big|\nabla_{r,g}g\big|^pr^{p-N}\dv.
\end{align*}
Now using \eqref{anoth-crit-cor} we get 
\begin{align*}
    J(f)&\geq c_p \frac{(p-1)^p}{p^p}\int_{\m}\frac{|g(x)-g(R\frac{x}{r})|^p}{r^{N}\big|\ln\frac{R}{r}\big|^{p}}\dv\\&=c_p \frac{(p-1)^p}{p^p}\int_{\m}\frac{||x|^{\frac{N-p-\beta}{p}}f(x)-R^{\frac{N-p-\beta}{p}}f(R\frac{x}{|x|})|^p}{|x|^{N}\big|\ln\frac{R}{|x|}\big|^{p}}\dv,
\end{align*}
for any $R>0$. After simplifying it, we arrive at
\begin{align*}
    J(f)\geq c_p \frac{(p-1)^p}{p^p}\int_{\m}\frac{|f(x)-R^{\frac{N-p-\beta}{p}}f(R\frac{x}{|x|})|x|^{-\frac{N-p-\beta}{p}}|^p}{|x|^{p+\beta}\big|\ln\frac{R}{|x|}\big|^{p}}\dv=c_p \frac{(p-1)^p}{p^p}d_H(f,R)^p.
\end{align*}
Now taking the supremum over $R>0$, we deduce the desired result.

{\bf Proof of Theorem~\ref{stab-th-crt}:} Let $x=r\sigma$ with $r=\rho(x)$. For any complex-valued functions $f\in C_0^\infty(\ba\setminus\{o\})$, we define the function $g(s\sigma)=\big(\ln \frac{1}{r}\big)^{-\frac{p-1}{p}}f(r\sigma)$, where $s=\big(\ln \frac{1}{r}\big)^{-1}$. We have $r\rightarrow 0$ for $s\rightarrow 0$ and so we get $g(o)=\lim_{s\rightarrow 0} g(s\sigma)=0$ and as $f$ has compact support so $g\in C_0^\infty(\m\setminus\{o\})$. 
By using the polar coordinate decomposition, we obtain 
\begin{align*}
&C(f):=\int_{\ba}\frac{|\nabla_{r,g} f|^p}{r^{N-p}}\dv - \bigg(\frac{p-1}{p}\bigg)^p\int_{\ba}\frac{|f|^p}{r^{N}(\ln \frac{1}{r})^p}\dv\\&=\int_{\sn}\int_0^1\bigg[\bigg|\frac{\partial f}{\partial r}\bigg|^p r^{p-1}-\bigg(\frac{p-1}{p}\bigg)^p\frac{|f|^p }{r(\ln \frac{1}{r})^p}\bigg]J(r,\sigma)\dr\dsn\\&=\bigg[\int_{\sn}\int_0^1\bigg|\frac{\partial g(s\sigma)}{\partial s} \frac{\partial s}{\partial r}\bigg(\ln \frac{1}{r}\bigg)^{\frac{p-1}{p}} \\&- \bigg(\frac{p-1}{p}\bigg)g(s\sigma) r^{-1}\bigg(\ln \frac{1}{r}\bigg)^{-\frac{1}{p}}\bigg|^pr^{p-1}J(r,\sigma)\dr\dsn\bigg]\\&-\bigg(\frac{p-1}{p}\bigg)^p\int_{\sn}\int_0^1\frac{|g(s\sigma)|^p }{r(\ln \frac{1}{r})}J(r,\sigma)\dr\dsn.
\end{align*}

Now applying Lemma~\ref{cpl}, we have 
\begin{align*}
    &\bigg|\bigg(\frac{p-1}{p}\bigg)g(s\sigma) r^{-1}\bigg(\ln \frac{1}{r}\bigg)^{-\frac{1}{p}}-\frac{\partial g(s\sigma)}{\partial s} \frac{\partial s}{\partial r}\bigg(\ln \frac{1}{r}\bigg)^{\frac{p-1}{p}}\bigg|^p\\&\geq \bigg(\frac{p-1}{p}\bigg)^p|g(s\sigma)|^pr^{-p}\bigg(\ln \frac{1}{r}\bigg)^{-1}-p\bigg(\frac{p-1}{p}\bigg)^{p-1}|g(s\sigma)|^{p-2}{\rm Re\;}\overline{g(s\sigma)}\frac{\partial g(s\sigma)}{\partial s} \frac{\partial s}{\partial r}r^{1-p}\\&+c_p\bigg|\frac{\partial g(s\sigma)}{\partial s}\bigg|^p\bigg|\frac{\partial s}{\partial r}\bigg|^p\bigg(\ln \frac{1}{r}\bigg)^{p-1}.
\end{align*}
Next, using integration by parts and the increasing property of $J(r,\sigma)$ with respect to $r$, we have the following
\begin{align*}
&-p\bigg(\frac{p-1}{p}\bigg)^{p-1}\int_{\sn}\int_0^1|g(s\sigma)|^{p-2}{\rm Re\;}\overline{g(s\sigma)}\frac{\partial g(s\sigma)}{\partial s} \frac{\partial s}{\partial r}r^{1-p}r^{p-1} J(r,\sigma)\dr\dsn\\&=-\bigg(\frac{p-1}{p}\bigg)^{p-1}\int_{\sn}\int_0^1\frac{\partial |g(s\sigma)|^p}{\partial r} J(r,\sigma)\dr\dsn\\&=\bigg(\frac{p-1}{p}\bigg)^{p-1}\int_{\sn}\int_{0}^1|g(s\sigma)|^pJ_r(r,\sigma)\dr\dsn \geq 0.
\end{align*}
In the above, we used 
\begin{align*}
    \left[\int_{\sn}|g(s\sigma)|^p J(r,\sigma)\dsn\right]_{r=0}^{r=1}=\left[\int_{\sn}|\big(\ln \frac{1}{r}\big)^{-\frac{p-1}{p}}f(r\sigma)|^p J(r,\sigma)\dsn\right]_{r=0}^{r=1}=0.
\end{align*}
The above follows because $g(s\sigma)=\big(\ln \frac{1}{r}\big)^{-\frac{p-1}{p}}f(r\sigma)$ and $f\in C_0^\infty(\ba\setminus\{o\})$. Therefore, we deduce the remainder term as below
\begin{align}\label{ref-step}
   \nonumber C(f)&\geq c_p\int_{\sn}\int_0^
    1 \bigg|\frac{\partial g(s\sigma)}{\partial s}\bigg|^p\bigg|\frac{\partial s}{\partial r}\bigg|^p\bigg(\ln \frac{1}{r}\bigg)^{p-1}r^{p-1}J(r,\sigma)\dr\dsn\\&= c_p\int_{\sn}\int_0^
    \infty \bigg|\frac{\partial g(s\sigma)}{\partial s}\bigg|^ps^{p-1}J(e^{-s^{-1}}, \sigma)\ds\dsn.
\end{align}
Arguing as in \eqref{steps-ref}, and using the fact that for some fixed $\sigma \in \sn$ the function $J_s(e^{-s^{-1}}, \sigma) \geq 0$ with respect to $s$, we obtain
\begin{align*}
    \int_{0}^R \frac{|g(s\sigma)-g(R\sigma)|^p}{s^{N}\big|\ln\frac{R}{s}\big|^{p}}s^{N-1}J(e^{-s^{-1}},\sigma)\ds\leq \left(\frac{p}{p-1}\right)^p\int_0^R
     \bigg|\frac{\partial g(s\sigma)}{\partial s}\bigg|^ps^{p-1}J(e^{-s^{-1}}, \sigma)\ds,
\end{align*}
and 
\begin{align*}
    \int_{R}^\infty \frac{|g(s\sigma)-g(R\sigma)|^p}{s^{N}\big|\ln\frac{R}{s}\big|^{p}}s^{N-1}J(e^{-s^{-1}},\sigma)\ds\leq \left(\frac{p}{p-1}\right)^p\int_R^\infty
     \bigg|\frac{\partial g(s\sigma)}{\partial s}\bigg|^ps^{p-1}J(e^{-s^{-1}}, \sigma)\ds,
\end{align*}
where $R > 0$ is any real number. Combining both of the above estimates, integrating over $\sn$, and substituting into \eqref{ref-step}, we arrive at
\begin{align*}
    C(f)\geq c_p \frac{(p-1)^p}{p^p}\int_{\sn}\int_0^\infty\frac{|g(s\sigma)-g(R\sigma)|^p}{s^{N}\big|\ln\frac{R}{s}\big|^{p}}s^{N-1}J(e^{-s^{-1}},\sigma)\ds\dsn.
\end{align*}

Now, in the above, applying the change of variable $r = e^{-s^{-1}}$ and using the identity  $g(s\sigma)=\big(\ln \frac{1}{r}\big)^{-\frac{p-1}{p}}f(r\sigma)=s^{\frac{p-1}{p}}f(e^{-s^{-1}}\sigma)$, we deduce
\begin{align*}
    &C(f)\geq c_p \frac{(p-1)^p}{p^p}\int_{\sn}\int_0^1\frac{|\big(\ln \frac{1}{r}\big)^{-\frac{p-1}{p}}f(r\sigma)-R^{\frac{p-1}{p}}f(e^{-R^{-1}}\sigma)|^p}{\big(\ln \frac{1}{r}\big)\big|\ln \big(R\ln \frac{1}{r}\big)\big|^{p}}r^{-1}J(r,\sigma)\dr\dsn\\&=c_p \frac{(p-1)^p}{p^p}\int_{\sn}\int_0^1\frac{|f(r\sigma)-\big(\ln \frac{1}{r}\big)^{\frac{p-1}{p}}R^{\frac{p-1}{p}}f(e^{-R^{-1}}\sigma)|^p}{r^N\big(\ln \frac{1}{r}\big)^{p}\big|\ln \big(R\ln \frac{1}{r}\big)\big|^{p}}r^{N-1}J(r,\sigma)\dr\dsn\\&=c_p \frac{(p-1)^p}{p^p}\int_{\ba}\frac{|f(x)-\big(\ln \frac{1}{|x|}\big)^{\frac{p-1}{p}}R^{\frac{p-1}{p}}f(e^{-R^{-1}}\sigma)|^p}{|x|^N\big(\ln \frac{1}{|x|}\big)^{p}\big|\ln \big(R\ln \frac{1}{|x|}\big)\big|^{p}}\dv\\&=c_p \frac{(p-1)^p}{p^p}d_C(f,R)^p,
\end{align*}
for any $R>0$. Finally, taking the supremum over $R>0$, we complete the proof.

\medspace

\section{Proofs of Theorem~\ref{new-hardy-th}, Corollary~\ref{ckn-th} and Corollary~\ref{rel-cor}}\label{pfs-4}

{\bf Proof of Theorem~\ref{new-hardy-th}:}   Denote $r=\rho(x)$. Let us start with
    \begin{align*}
       (b-1)c &\int_{\ba} \frac{|f|^p}{r^a(1-r^c)^b}\dv \\& =\int_{\ba}\frac{|f|^p}{r^{a+(b-1)c-1}}\frac{{\rm d}}{{\rm d}r}\bigg(\frac{1}{(r^{-c}-1)^{b-1}}\bigg)\dv \\&=\int_{\sn}\int_{0}^{1} \frac{|F(r\sigma)|^p}{r^{a+(b-1)c-N}}\frac{{\rm d}}{{\rm d}r}\bigg(\frac{1}{(r^{-c}-1)^{b-1}}\bigg) J(r,\sigma) \dr\dsn\\&=-p{\rm Re}\int_{\sn}\int_{0}^{1} \frac{|F|^{p-2}F \overline{F_r}}{r^{a+(b-1)c-N}}\frac{1}{(r^{-c}-1)^{b-1}} J(r,\sigma) \dr\dsn\\&-\int_{\sn}\int_{0}^{1} \frac{|F|^p}{r^{a+(b-1)c-N}}\frac{1}{(r^{-c}-1)^{b-1}} J_r(r,\sigma) \dr\dsn\\&-\big(N-a-(b-1)c\big)\int_{\sn}\int_{0}^{1} \frac{|F|^p}{r^{a+(b-1)c-N+1}}\frac{1}{(r^{-c}-1)^{b-1}} J(r,\sigma) \dr\dsn\\&=-p{\rm Re}\int_{\ba} \frac{|f|^{p-2}\, f \, \overline{\frac{\partial f}{\partial r}}}{r^{a+(b-1)c-1}}\frac{1}{(r^{-c}-1)^{b-1}} \dv\\&-\int_{\ba}\frac{|f|^p}{r^{a+(b-1)c-1}}\frac{1}{(r^{-c}-1)^{b-1}} \frac{J_r(r,\sigma)}{J(r,\sigma)} \dv\\&-\big(N-a-(b-1)c\big)\int_{\ba} \frac{|f|^p}{r^{a+(b-1)c}}\frac{1}{(r^{-c}-1)^{b-1}} \dv\\&=-p{\rm Re}\int_{\ba} \frac{|f|^{p-2}\, f \, \overline{\frac{\partial f}{\partial r}}}{r^{a-1}}\frac{1}{(1-r^{c})^{b-1}} \dv\\&-\int_{\ba}\frac{|f|^p}{r^{a-1}}\frac{1}{(1-r^{c})^{b-1}} \frac{J_r(r,\sigma)}{J(r,\sigma)} \dv\\&-\big(N-a-(b-1)c\big)\int_{\ba} \frac{|f|^p}{r^{a}}\frac{1}{(1-r^{c})^{b-1}} \dv\\&\leq p\bigg(\int_{\ba} \frac{|\frac{\partial f}{\partial r}|^p}{r^{a-p}(1-r^{c})^{b-p}} \dv\bigg)^{\frac{1}{p}}\bigg(\int_{\ba} \frac{|f|^{p}}{r^{a}(1-r^{c})^{b}} \dv\bigg)^{1-\frac{1}{p}}\\&-\big(N-a-(b-1)c\big)\int_{\ba} \frac{|f|^p}{r^{a}(1-r^{c})^{b-1}} \dv.
    \end{align*}
    Also, in the above, we used 
    \begin{align*}
        \left[\int_{\sn} \frac{|F(r\sigma)|^p}{r^{a+(b-1)c-N}}\frac{1}{(r^{-c}-1)^{b-1}} J(r,\sigma) \dsn\right]_{r=0}^{r=1}=0.
    \end{align*}
    It can be verified at $r=0$ using the condition $a<N$ and at $r=1$ using the compact support of $f$. 
   Therefore,
    \begin{align*}
        \frac{(b-1)c}{p}\bigg(\int_{\ba} \frac{|f|^{p}}{r^{a}(1-r^{c})^{b}} \dv\bigg)^{\frac{1}{p}}&\leq \bigg(\int_{\ba} \frac{|\frac{\partial f}{\partial r}|^p}{r^{a-p}(1-r^{c})^{b-p}} \dv\bigg)^{\frac{1}{p}}\\&- \frac{\big(N-a-(b-1)c\big)}{p}\frac{\int_{\ba} \frac{|f|^p}{r^{a}(1-r^{c})^{b-1}} \dv}{\bigg(\int_{\ba} \frac{|f|^{p}}{r^{a}(1-r^{c})^{b}} \dv\bigg)^{\frac{p-1}{p}}}.
    \end{align*}   
   Omitting the latter negative term we obtain (for $f\not \equiv 0 $) the following inequality 
    \begin{align*}
        \bigg(\frac{(b-1)}{p}c\bigg)^p\int_{\ba} \frac{|f|^{p}}{r^{a}(1-r^{c})^{b}} \dv\leq \int_{\ba} \frac{|\frac{\partial f}{\partial r}|^p}{r^{a-p}(1-r^{c})^{b-p}} \dv.
    \end{align*}
    Now we verify the optimality of the constant. Take $t>\frac{b-1}{p}$, and small positive $\delta$ such that $2\delta\in (0,1)$. Let us define the following radially symmetric function
    \begin{align*}
        f_t(r)=\chi(r) (1-r^c)^t,
    \end{align*}
    where cut-off function $\chi:[0,\infty)\rightarrow \mathbb{R}$ satisfies the following properties:
\begin{itemize}
    \item[1.] $\chi(r)\in [0,1]$ for all $r\in [0,\infty)$ and $\chi$ is smooth and compactly supported;
    \item[2.] 
\begin{equation*}
\chi(r)=
\begin{dcases}
1, & 1-\delta\leq r\leq 1; \\
0, & 0\leq r< 1-2\delta. \\
\end{dcases}
\end{equation*}
\end{itemize}
It is easy to verify that $f_t \in L^{p}(\ba)$, and by using the condition $p \leq b$ together with the support of the cutoff function, we can demonstrate that $(f_t)^{\prime} \in L^{p}(\ba)$. This implies that $f_t \in W^{1,p}_{0}(\ba)$, and these will serve as test functions due to the denseness of smooth functions in these Sobolev spaces. By using (increasing) monotonicity of $J(r,\sigma)$ and \eqref{imp-den} we deduce
\begin{align*}
    \int_{\ba} &\frac{|f_t|^{p}}{r^{a}(1-r^{c})^{b}} \dv  = \int_{\sn}\int_{1-2\delta}^{1} \frac{|\chi(r) (1-r^c)^t|^{p}}{r^{a}(1-r^{c})^{b}} r^{N-1}J(r,\sigma)\dr\dsn \\& \geq \int_{\sn}\int_{1-\delta}^{1} r^{N-1-a}(1-r^c)^{tp-b}J(r,\sigma)\dr\dsn\\&\geq |\sn|\int_{1-\delta}^{1} r^{N-1-a}(1-r^c)^{tp-b}\dr\\&=|\sn|\bigg(-\frac{1}{c}\bigg)\int_{1-\delta}^{1} r^{N-c-a}(1-r^c)^{tp-b}\,{\rm d}(1-r^c)\\&\geq |\sn|\bigg(-\frac{1}{c}\bigg)\max\{1,(1-\delta)^{N-c-a}\}\int_{1-\delta}^{1} (1-r^c)^{tp-b}\,{\rm d}(1-r^c)\\&=|\sn|\bigg(-\frac{1}{c(tp-b+1)}\bigg)\max\{1,(1-\delta)^{N-c-a}\}\bigg[(1-r^c)^{tp-b+1}\bigg]_{1-\delta}^1\\&=|\sn|\max\{1,(1-\delta)^{N-c-a}\}\bigg[\frac{(1-(1-\delta)^c)^{tp-b+1}}{c(tp-b+1)}\bigg],
    \end{align*}
which tends to $\infty$ as $t\rightarrow \frac{b-1}{p}$. Therefore,
\begin{align*}
  \lim_{t\rightarrow \frac{b-1}{p}}  \int_{\ba} \frac{|f_t|^{p}}{r^{a}(1-r^{c})^{b}} \dv =\infty.
\end{align*}
On the other hand, we have
\begin{align*}
  \frac{\partial f_t}{\partial r}= \chi^\prime(r)(1-r^c)^t-tc\chi(r)(1-r^c)^{t-1}r^{c-1} .
\end{align*}
Therefore, we establish 
\begin{align*}
   &\int_{\ba} \frac{|\chi^\prime(r)(1-r^c)^t|^p}{r^{a-p}(1-r^{c})^{b-p}}\dv\\&=\int_{\sn}\int_{1-2\delta}^{1-\delta}|\chi^\prime(r)|^pr^{N-a+p-1}(1-r^c)^{tp-b+p}J(r,\sigma)\dr\dsn=O(1),
\end{align*}
and 
\begin{align*}
   &(tc)^p\int_{\ba} \frac{|\chi(r)(1-r^c)^{t-1}r^{c-1}|^p}{r^{a-p}(1-r^{c})^{b-p}}\dv\\&=(tc)^p\int_{\ba} r^{pc}\frac{|\chi(r)(1-r^c)^{t}|^p}{r^{a}(1-r^{c})^{b}}\dv\\& \leq (tc)^p \int_{\ba} \frac{|f_t|^{p}}{r^{a}(1-r^{c})^{b}} \dv, 
\end{align*}
for $t\rightarrow \frac{b-1}{p}$. That is, 
\begin{align*}
    \lim_{t\rightarrow \frac{b-1}{p}}\frac{\int_{\ba} \frac{|\frac{\partial f_t}{\partial r}|^p}{r^{a-p}(1-r^{c})^{b-p}} \dv}{\int_{\ba} \frac{|f_t|^{p}}{r^{a}(1-r^{c})^{b}} \dv}=\lim_{t\rightarrow \frac{b-1}{p}} (tc)^p= \bigg(\frac{(b-1)}{p}c\bigg)^p.
\end{align*}

{\bf Proof of Corollary~\ref{ckn-th}:} For the case $\delta=0$, the inequality \eqref{ckn} becomes trivial as we have $\eta=q$ and $\gamma=\beta$. When $\delta=1$, we have $\eta=p$ and $\gamma=\alpha-1$. By using these facts and \eqref{new-wg-hardy}, we establish \eqref{ckn}. Now the case $\delta\in (0,1)\cap \big[\frac{\eta-q}{\eta}, \frac{p}{\eta}\big]$ is remained. One can simply write
\begin{align*}
    ||w^\gamma f||_{L^\eta(\ba)}=\int_{\ba}\frac{|f(x)|^{\delta\eta}}{(w(x))^{\delta\eta(1-\alpha)}}\cdot \frac{|f(x)|^{(1-\delta)\eta}}{(w(x))^{-\beta\eta(1-\delta)}} \dv.
\end{align*}
Then by using the H\"older inequality for the exponent $\frac{\delta \eta}{p}+\frac{(1-\delta)\eta}{q}=1$ we deduce
\begin{align*}
   ||w^\gamma f||_{L^\eta(\ba)}\leq ||w^{\alpha-1} f||_{L^p(\ba)}^\delta||w^{\beta} f||_{L^q(\ba)}^{1-\delta}.
\end{align*}
Hence, applying \eqref{new-wg-hardy} on the first term, we obtain \eqref{ckn}.

{\bf Proof of Corollary~\ref{rel-cor}:} With $r=\rho(x)$ applying \eqref{rell-typ} for the function $\frac{\partial f}{\partial r}$, we have
\begin{align*}
\bigg(\frac{N(p-1)+\beta}{p}\bigg)^p\int_{\ba}\frac{|\nabla_{r,g} f|^p}{r^{p+\beta}}\dv\leq \int_{\ba}\frac{|\Delta_{r,g} f|^p}{r^\beta}\dv.
\end{align*}
Now by using \eqref{new-wg-hardy} on the left-hand side, we get
\begin{align*}
        \bigg(\frac{(p-1)}{p}c\bigg)^p\int_{\ba} \frac{|f|^{p}}{r^{2p+\beta}(1-r^{c})^{p}} \dv\leq \int_{\ba} \frac{|\nabla_{r,g} f|^p}{r^{p+\beta}} \dv.
\end{align*}
Combining these two completes the proof of the first part. Moreover, the constant is optimal here as both the constants in the respective inequalities are sharp.

In particular choosing $c=\frac{N-2p-\beta}{p-1}$, and for $0<r<1$, using the inequality
\begin{align*}
    1\leq \big(1-r^{\frac{N-2p-\beta}{p-1}}\big)^{-1},
\end{align*}
finishes the proof of the second part. The optimality follows from the known optimal constant in this setting \cite[Theorem~4.3]{vhn}.

\medspace

\section{Proof of Theorem~\ref{higher-order}}\label{pfs-5}
First, notice that for $k=1$, taking $b=p$ and $a=p+\beta$ in Theorem~\ref{new-hardy-th}, we deduce the result. Now observe that the case $k=2$ follows from Corollary~\ref{rel-cor}. So, we only present the cases when $k=3$ and $k=4$ here. Other cases will simply follow by considering induction in the odd and even cases separately.

{\bf The case $k=3$:} Assume $r=\rho(x)$. Let us start with applying \eqref{new-wg-hardy}, for $a-p=\beta$ and $b=p$ for the function $\Delta_{r,g} f$, and we have
    \begin{align*}
     \int_{\ba} \frac{|\nabla_{r,g}^3 f|^p}{r^{\beta}} \dv  & \geq \bigg(\frac{(p-1)}{p}c\bigg)^p\int_{\ba} \frac{|\Delta_{r,g}f|^{p}}{r^{p+\beta}(1-r^{c})^{p}} \dv\\& \geq \bigg(\frac{(p-1)}{p}c\bigg)^p\int_{\ba} \frac{|\Delta_{r,g}f|^{p}}{r^{p+\beta}} \dv.
    \end{align*}
    By applying \eqref{rel-cor-eqn}, we get
    \begin{align*}
     \int_{\ba} \frac{|\nabla_{r,g}^3 f|^p}{r^{\beta}} \dv  & \geq \bigg(\frac{(p-1)}{p}c\bigg)^p\int_{\ba} \frac{|\Delta_{r,g}f|^{p}}{r^{p+\beta}} \dv \\& \geq \bigg(\frac{N(p-1)+\beta+p}{p}\bigg)^p\bigg(\frac{(p-1)}{p}c\bigg)^{2p}\int_{\ba} \frac{|f|^{p}}{r^{3p+\beta}(1-r^{c})^{p}} \dv.
    \end{align*}

{\bf The case $k=4$:} Choose $r=\rho(x)$. By using Corollary~\ref{rel-cor} twice for the function $\Delta_{r,g} f$ and then for $f$, we obtain
    \begin{align*}
     &\int_{\ba} \frac{|\nabla_{r,g}^4 f|^p}{r^{\beta}} \dv  \\& \geq \bigg(\frac{N(p-1)+\beta}{p}\bigg)^p\bigg(\frac{(p-1)}{p}c\bigg)^p\int_{\ba} \frac{|\Delta_{r,g} f|^{p}}{r^{2p+\beta}(1-r^{c})^{p}} \dv\\& \geq \bigg(\frac{N(p-1)+\beta}{p}\bigg)^p\bigg(\frac{(p-1)}{p}c\bigg)^p\int_{\ba} \frac{|\Delta_{r,g}f|^{p}}{r^{2p+\beta}} \dv \\& \geq \bigg(\frac{N(p-1)+\beta+2p}{p}\bigg)^p\bigg(\frac{N(p-1)+\beta}{p}\bigg)^p\bigg(\frac{(p-1)}{p}c\bigg)^{2p}\int_{\ba} \frac{|f|^{p}}{r^{4p+\beta}(1-r^{c})^{p}} \dv.
    \end{align*}

\medspace

\section{Open problem}\label{sec-open}
In Theorem~\ref{stab-th} and Theorem~\ref{stab-th-crt}, we discuss the stability inequalities for the $L^p$-Hardy inequality in both the subcritical and critical cases. It would be interesting to study the optimality of the constant that occurs in front of the stable metric parameter. Here we formulate our question for the subcritical case (a similar question can be asked for the critical case): Let $f\in C_0^\infty (\m\setminus \{o\})$ and define
\begin{align*}
		E(f):=\frac{\int_{\m}\frac{|\partial_\rho f(x)|^p}{\rho^\beta(x)}\dv- \big(\frac{N-p-\beta}{p}\big)^p\int_{\m}\frac{|f(x)|^p}{\rho^{p+\beta}(x)}\dv}{\sup_{R>0} d_H(f,R)^p}.
	\end{align*}
Theorem~\ref{stab-th} gives
$$ E(f)\geq c_p\bigg(\frac{p-1}{p}\bigg)^p.$$
Does $\inf_{f\in C_0^\infty (\m\setminus \{o\})}E(f)= c_p\big(\frac{p-1}{p}\big)^p$ hold? This question appears to remain open even in the Euclidean setting, not only in the context of Cartan--Hadamard manifolds.

\medspace

\section*{Acknowledgments}
 The authors are grateful to the anonymous referee for their insightful comments, which have substantially improved the clarity and presentation of the manuscript, especially Theorem~\ref{stab-th-crt}. This research was funded by the Committee of Science of the Ministry of Science and Higher Education of  Kazakhstan (Grant No. AP23488549). This project was completed when the authors met at the Department of Mathematics at Nazarbayev University (NU) in the Autumn of 2023. The first author is partially supported by the National Theoretical Science Research Center Operational Plan (Project number: 112L104040) and the MUR-PRIN project No. 20227HX33Z ``Pattern formation in nonlinear phenomena'' granted by the European Union - Next Generation EU. The second and third authors were supported by the NU project 20122022FD4105. The third author is also supported by the Committee of Science of the Ministry of Science and Higher Education of Kazakhstan (Grant No. AP23490970).   

\medspace

\end{document}